\documentclass[a4paper,11pt]{article}

\usepackage{color,url}
\usepackage{enumitem}
\usepackage{graphicx}

\usepackage{psfrag}
\usepackage{epic}
\usepackage{amssymb}
\usepackage{afterpage}
\usepackage{epsfig}
\usepackage{pstricks}
\usepackage{stmaryrd,url}
\usepackage{tikz,amsmath}
\usepackage{diagbox}
\usepackage{framed,multirow}
\usepackage{epstopdf}

\usepackage[top=1.0in,bottom=1.0in,left=1.0in,right=1.0in]{geometry}
\usepackage{hyperref}

\newcommand{\ba}{\begin{array}}
	\newcommand{\ea}{\end{array}}
\newcommand{\be}{\begin{equation}}
\newcommand{\ee}{\end{equation}}
\newcommand{\ben}{\begin{equation*}}
\newcommand{\een}{\end{equation*}}
\newcommand{\bd}{\begin{displaymath}}
\newcommand{\ed}{\end{displaymath}}
\newcommand{\bi}{\begin{itemize}}
	\newcommand{\ei}{\end{itemize}}
\newcommand{\bn}{\begin{enumerate}}
	\newcommand{\en}{\end{enumerate}}
\newcommand{\pa}{\partial}
\newcommand{\f}{\frac}
\newcommand{\ci}{\cite}

\title{Energy-based discontinuous Galerkin difference methods for second-order wave equations\thanks{This work  
	was supported by NSF Grants DMS-1913076 and DMS-2012296 and completed while the third author was in residence at the Institute for Computational and
	Experimental Mathematics. Any opinions, findings, and conclusions or recommendations expressed in this material are those of the authors and do not necessarily reflect the views of the National Science Foundation}}

\begin{document}
\author{Lu Zhang\thanks{Department of Applied Physics and Applied Mathematics, Columbia University, New York, NY. \href{mailto:lz2784@columbia.edu}{Email: lz2784@columbia.edu}}\and Daniel  Appel\"{o}\thanks{Department of Computational Mathematics, Science \& Engineering and Department of Mathematics,
		Michigan State University, East Lansing, MI. \href{mailto:appeloda@msu.edu}{Email: appeloda@msu.edu}}\and Thomas Hagstrom\thanks{Department of Mathematics, Southern Methodist University, Dallas, TX. \href{mailto:thagstrom@smu.edu}{Email: thagstrom@smu.edu}}}
\maketitle

\begin{abstract}
We combine the newly-constructed Galerkin difference basis with the energy-based discontinuous Galerkin method for wave equations in second order form. The approximation properties of the resulting method are excellent and the allowable time steps are large compared to traditional discontinuous Galerkin methods. The one drawback of the combined approach is the cost of inversion of the local mass matrix. We demonstrate that for constant coefficient problems on Cartesian meshes this bottleneck can be removed by the use of a modified Galerkin difference basis. For variable coefficients or non-Cartesian meshes this technique is not possible and we instead use the preconditioned conjugate gradient method to iteratively invert the mass matrices. With a careful choice of preconditioner we can demonstrate optimal complexity, albeit with a larger constant.    
\end{abstract}


\textbf{AMS subject }: 65M60,  65M06 

\section{Introduction}

Discontinuous Galerkin methods have become a method of choice for solving first order hyperbolic equations in Friedrichs form \cite{hesthaven2007nodal}. They possess many desirable properties such as arbitrary order, robustness, geometric flexibility and explicit time evolution.  Analogous methods for second order hyperbolic equations are less well established, despite the fact that many governing equations arising in physics are in second order form. Even though it is often possible to rewrite second order hyperbolic equations in first order form, the first order formulation has some drawbacks. It almost always needs more
boundary conditions and it is only equivalent to the original second order equation for constrained data. In addition, not all second order hyperbolic equations can be rewritten as a first order system which is in Friedrichs form. In \cite{DATH_UP} Appel{\"o} and Hagstrom  proposed an energy-based DG method for second order wave equations. The idea in \cite{DATH_UP} is to introduce a new variable $v = \frac{\partial u}{\partial t}$ to transfer second order hyperbolic equations to first order systems in time only, and then seek approximations which satisfy a discrete energy equality. The method features a direct, mesh-independent approach to defining interelement fluxes. Both energy-conserving and upwind discretizations have been devised and their extension to elastic, advective and semi-linear wave equations can be found in \cite{appelo2018energy,appelo2020energy,zhang2019energy}.

Galerkin difference methods were introduced by Banks and Hagstrom \cite{Gdiff} to solve hyperbolic initial-boundary value problems. The idea is to use a Galerkin construction to derive energy stable finite difference methods. The basis functions are Lagrange functions associated with continuous piecewise polynomial approximation on a computational grid. Salient features of these methods are: the discrete approximations are uniform from grid-point to grid-point when away from domain boundaries; no new degrees of freedom are added when the approximation order increases; they do not require significantly smaller time steps as the order increases. Comparing them to summation-by-parts (SBP) difference schemes \ci{SBPrev}, they have the advantage of being directly constructable at arbitrary order and can be seamlessly interfaced with standard schemes on unstructured grids
\ci{GDgeo}. The relative disadvantage is that Galerkin difference operators typically require twice as many flops as SBP operators of the same order. 

In this paper, we combine the energy based discontinuous Galerkin methods with the Galerkin difference methods to solve the second order wave equation. The corresponding mass and stiffness matrices are banded matrices because of the properties of the Galerkin difference basis functions. For high dimensional problems on structured grids, the Galerkin basis functions are the tensor product of the Galerkin basis functions in one dimension. Then the mass matrix is a Kronecker product of the mass matrices in each dimension and the stiffness matrix is a summation of the Kronecker product of the mass and stiffness matrices in each dimension. This fact indicates that the inversion of the mass matrix can be computed with a computational cost which is linear with respect to the total degrees of freedom, while the computational cost grows rapidly for the inversion of the stiffness matrix. A result of this paper is the application of the simultaneous diagonalization technique from \cite{lynch1964direct} to derive a new class of Galerkin difference basis functions to reduce the computational cost of the inversion of the stiffness matrix.

We note that the literature on high order methods for wave equations in second order form is extensive. We will not try to review all methods here, but rather mention a few that are representative of the state of the art. In the class of discontinuous Galerkin methods it is worth mentioning the symmetric interior penalty (SIPG) method of \ci{GSSwave}, 
the local discontinuous Galerkin (LDG) method of \ci{ChouShuXing2014} and the nonsymmetric interior penalty 
method of \ci{IPDG_Elastic}. Finite difference methods include those using SBP operators \cite{Mattsson2004,Mattsson2012} as well as those using upwind methods \cite{BANKS20125854}. 

The rest of the paper is organized as follows. In Section~\ref{preli}, we introduce the construction of the Galerkin difference basis functions and the energy based discontinuous Galerkin method for the second order scalar wave equation. We derive the new basis functions from the Galerkin difference basis functions and investigate the reduction in computational cost in Section~\ref{GD_cost}. In Section~\ref{dispersion}, we present the dispersion properties of the scheme. Section~\ref{spectral_radius} shows the spectral radii of the proposed scheme. Numerical experiments that illustrate optimal convergence in both $L^2$ and energy norms are given in Section~\ref{convergence}. In Section \ref{variable}, we apply the method to  a problem with variable sound wave speed. Here the new basis construction does not apply, and as an alternative we use preconditioned conjugate gradient iterations. Finally, our conclusions summarized in Section~\ref{conclusion}.

\section{Preliminaries}\label{preli}
We consider the wave equation in first order form in time and second order form in space 
\begin{eqnarray}
u_t &=& v, \label{eq:wave_1} \\
v_t &=& \nabla \cdot (c^2 \nabla u) + f, \ \ t>0, \ \ (x,y,z) \in \Omega, \label{eq:wave_2}
\end{eqnarray}
complemented with initial data and boundary conditions. 

In what follows we will consider energy based discontinuous Galerkin methods implemented on $d$-dimensional tensor product elements. Each element will be mapped to the $d$-dimensional unit cube which will be discretized by an equidistant Cartesian grid. Precisely,
in each dimension $k$ and for each element we consider a grid discretizing the unit cube
\[
x_{k,i} = i h_k, \ \ i = 0,\ldots,N_k, \ \ h_k = 1/N_k.   
\]
For simplicity here we will use the same number of intervals in each dimension, $N_k=N$. In the following $n$ denotes the number of elements in one dimension, $m$ denotes the number of elements in multiple dimensions and $p$ denotes the polynomial degree of the approximation. We will only consider $p$ odd, which corresponds to the {\em locally continuous} difference basis and will often use the integer $q=(p+1)/2$.

\subsection{The Galerkin Difference Basis}\label{GD_basis}
We now describe the Galerkin difference (GD) basis we will use. Note that the description here differs slightly different from the original description in \cite{Gdiff}, but the basis is identical.    

We first consider grid points well-separated from a boundary. The goal is to construct a polynomial basis of odd degree $p$ on an equidistant grid with grid spacing $h$. The generating basis function $\Phi_p(x)$ is centered around $x=0$ and the basis itself is simply the union of the translates of the generating basis function. Thus an element in the basis centered around $x_i = ih$ becomes, $\phi_{p,i} = \Phi_p(x-ih)$. 
\graphicspath{{graphs_GDbasis/}}
\begin{figure}[h]
\begin{center}
\includegraphics[width=0.27\textwidth]{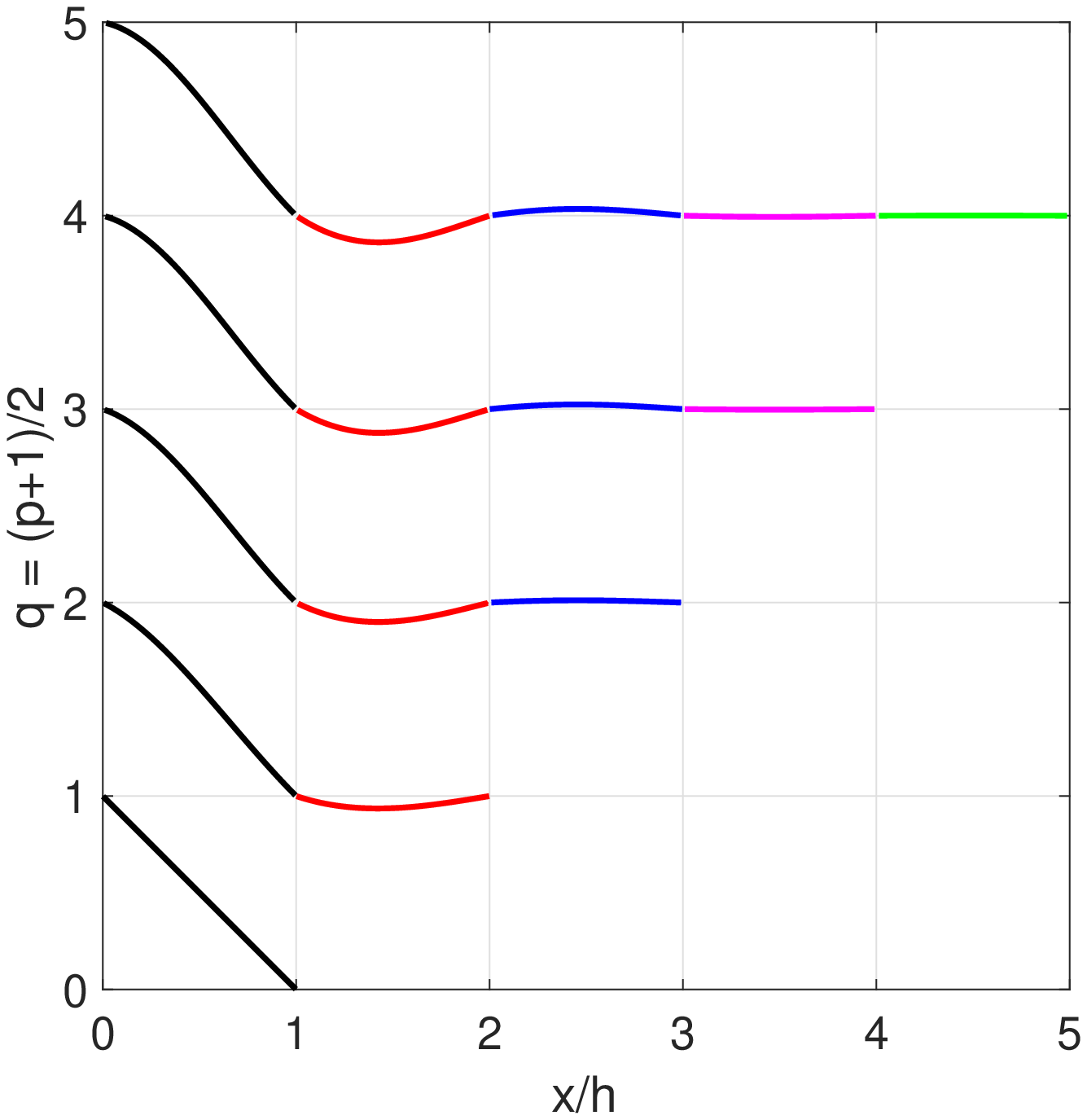}\ \ \ \
\includegraphics[width=0.21\textwidth]{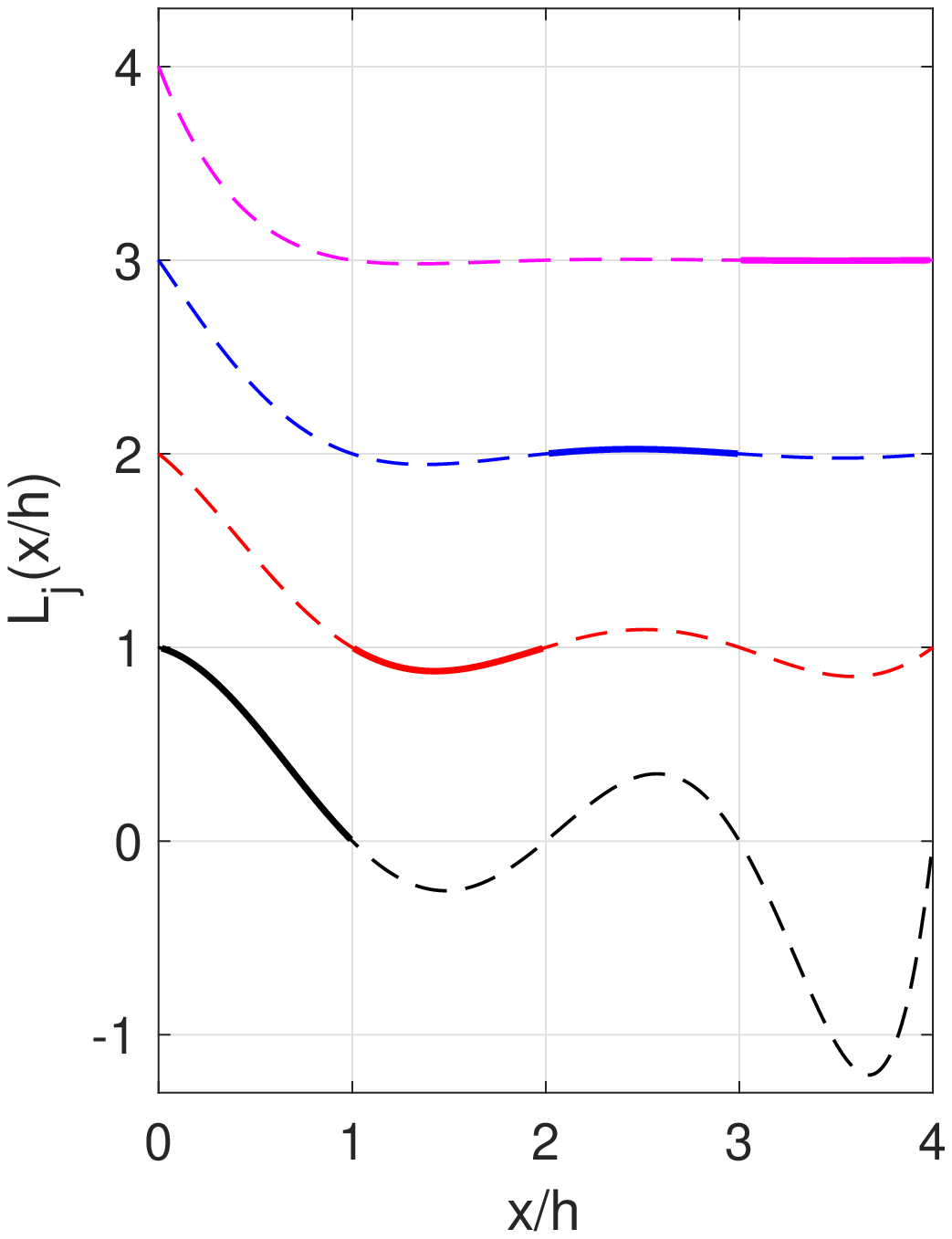} 
\includegraphics[width=0.41\textwidth]{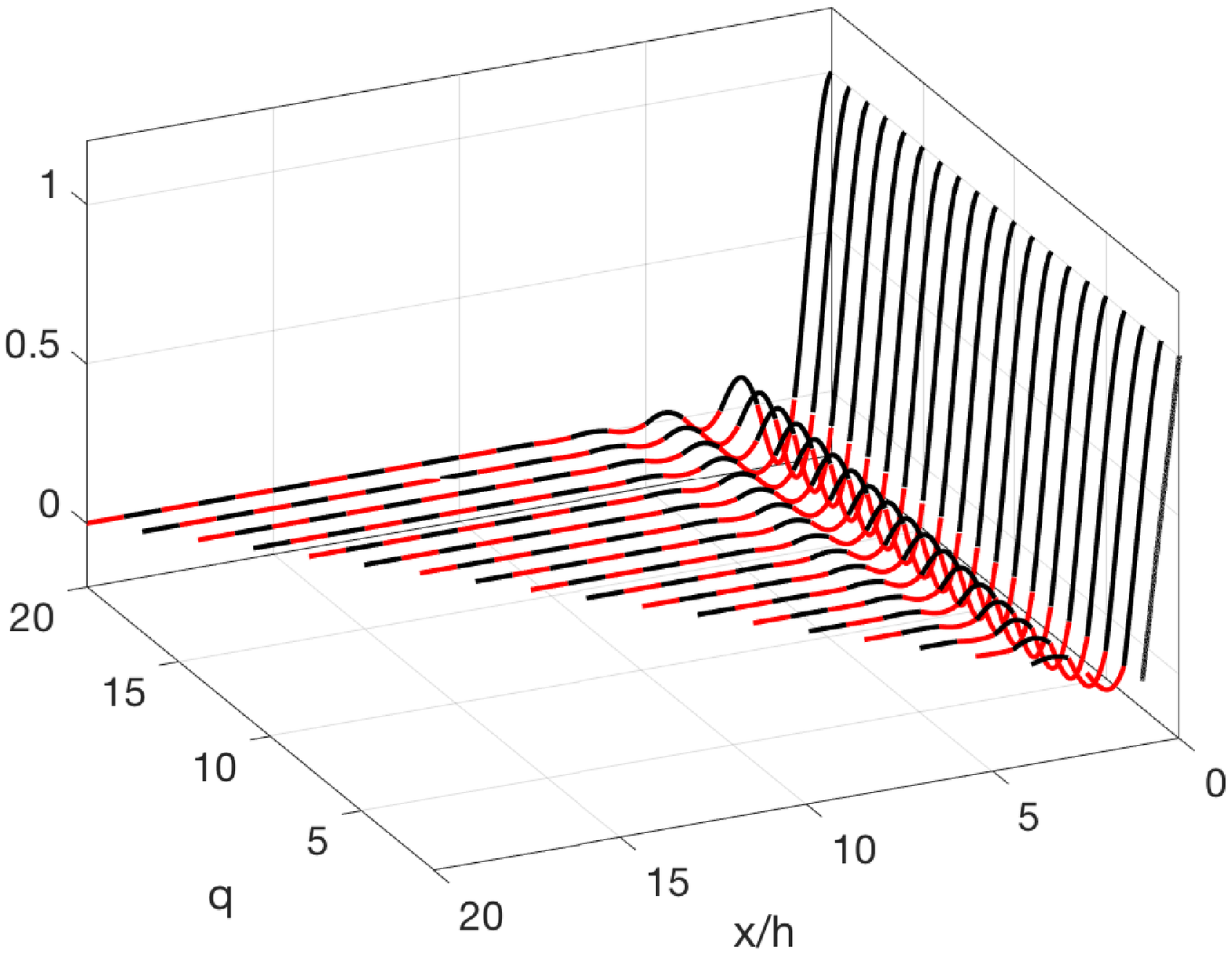}
\caption{To the left, the generating basis function $\Phi_p(x)$ for $p = 1,3,5,7,9$. Note that $\Phi_p(x) = \Phi_p(-x)$. Note that each function is vertically offset by $q-1$. In the middle, the ``right'' four degree 7 Lagrange polynomials that make up $\Phi_7(x)$, the part of $L_j$ that is used is in bold. To the right, a three dimensional representation of the generating basis function $\Phi_p(x)$ for $p$ up to 39. \label{fig:GD_generating}}
\end{center}
\end{figure}

The generating basis function $\Phi_p(x)$ is symmetric, $\Phi_p(x) = \Phi_p(-x)$, and has compact support on \mbox{$x \in [-qh,qh]$}, where $p = 2q-1$ (recall that $p$ is odd). In Figure \ref{fig:GD_generating} we display the non-zero part of $\Phi_p(x)$  for $p = 1,3,5,7,9,$ and $x>0$. In the lower left corner, where $p=1$, we recognize the classic finite element hat function and as $p$ increases we see that $\Phi_p(x)$ becomes increasingly similar to the {\em Cardinal Sinc} function.

An explicit formula for $\Phi_p(x)$ inside each of the $q$ positive intervals $[jh,(j+1)h)$ is obtained as follows. Let $L_j(x)$ is be the Lagrange interpolating polynomial on the grid $G_j = \{-jh,\ldots,(p-j)h\}$ with the property that $L_j(0) = 1$, then
\begin{eqnarray}\label{GD_basis_lagrange}
\Phi_p(x) = L_j(x), \ \ x \in [(q-j-1)h,(q-j)h)	. 
\end{eqnarray}

A continuous function $u(x,t)$ can then be approximated by a linear combination of basis functions with weights corresponding to nodal values
\begin{eqnarray}\label{GD_basis_nodalformula}
u(x,t) \approx \sum_{i=k-(q-1)}^{k+q} u_i \phi_{p,i}(x), \ \  x \in [kh,(k+1)h),
\end{eqnarray}
where $u_i = u(x_i,t)$.

\subsubsection{Modification Near Boundaries}\label{GD_basis_subsection_1}
Near boundaries the basis must be modified. In \cite{Gdiff}, three approaches for handling boundaries are described: ghost basis, extrapolation basis and use of modified equations. Here we will exclusively use the extrapolation basis, which we describe next.

The extrapolation procedure draws from the standard practice to use ghost points in finite difference methods. First, note that the $q-1$ additional ghost basis functions associated with the $q-1$ first grid points outside the computational domain are the only ghost basis with support inside the computational domain. In the ghost basis approach the degrees of freedom at the ghost points are retained as unknowns but in the extrapolation approach they are eliminated in favor of modifying the basis itself near the boundary. 

As the name suggests, the elimination is done by extrapolating the nodal values inside the computational domain to the ghost points. For example consider $p=3$. Then $q-1 = 1$ and one ghost point value, $u_{-1}$, must be determined. As the basis is fourth order accurate, the ghost point value is determined by fourth order accurate extrapolation,
\begin{eqnarray*}
u_{-1} &=& 4u_0 - 6u_1 + 4u_2 -u_3.
\end{eqnarray*}
 To understand how the modified basis is constructed, consider evaluating the approximation $u(x)$ inside the computational domain where the ghost basis has support. In this case this means \mbox{$x\in(x_0,x_1)$} and the approximation is
\begin{eqnarray*}
u(x) &=& u_{-1}\phi_{-1}(x) + u_{0}\phi_{0}(x) + u_{1}\phi_1(x) + u_2\phi_2(x).
\end{eqnarray*}
To obtain a value for $u_{-1}$ we use the extrapolation condition 
\begin{eqnarray*}
u(x) &=& ( 4u_0 - 6u_1 + 4u_2 -u_3)\phi_{-1}(x) + u_{0}\phi_{0}(x) + u_{1}\phi_1(x) + u_2\phi_2(x) \nonumber\\
&=& u_0 \left[ \phi_0(x)+4\phi_{-1}(x)\right] +u_1\left[ \phi_1(x)-6\phi_{-1}(x) \right] + u_2 \left[ \phi_2(x)+4\phi_{-1}(x) \right] + u_3 \left[-\phi_{-1}(x) \right] \\
& =& u_0 \left[ \phi_0(x)+4\phi_{-1}(x) \right] +u_1 \left[ \phi_1(x)-6\phi_{-1}(x) \right] + u_2 \left[ \phi_2(x)+4\phi_{-1}(x) \right] + u_3 \left[\phi_3(x)-\phi_{-1}(x ) \right] .\nonumber
\end{eqnarray*}
In the last step, we used the fact the support of $\phi_3$ vanishes in $(x_0,x_1)$. Thus the modified basis functions are 
\begin{equation*}
\tilde{\phi}_0 = \phi_0 + 4\phi_{-1}, \ \ \ \ 
\tilde{\phi}_1 = \phi_1 - 6\phi_{-1},   \ \ \ \ 
\tilde{\phi}_2 = \phi_2 + 4\phi_{-1},  \ \ \ \ 
\tilde{\phi}_3 = \phi_3 - \phi_{-1}.     
\end{equation*}
The extension to larger $p$ requires the basis to be modified in a wider band near the boundaries and the extrapolation is done at the order of accuracy that matches that of the interior scheme. The handling of the right boundary is analogous. 

\subsubsection{Extension to Higher Dimensions}\label{GD_basis_subsection_2}
The extension to higher dimensions simply amounts to using the tensor product approximation built off the one dimensional basis. For example, in two dimensions we have
\begin{eqnarray}
u(x,y,t) \approx \sum_{i=k_x-(q-1)}^{k_x+q}\sum_{j=k_y-(q-1)}^{k_y+q}u_{i,j} \phi_{p,i,j}(x,y), \ \ (x,y)\in[k_xh,(k_x+1)h)\times[k_yh,(k_y+1)h),
\end{eqnarray}
where $u_{i,j} = u(x_i,y_j,t)$ and 
\begin{equation}
\phi_{p,i,j}(x,y) = \left\{
\begin{array}{lll}
 \phi_{p,i}(x)\phi_{p,j}(y), & p <   i < N-p, & p < j < N - p, \\
 {\phi}_{p,i}(x)\tilde{\phi}_{p,j}(y), & p<i<N-p, & 0 \le j \le N,\\
 \tilde{\phi}_{p,i}(x)\phi_{p,j}(y), &  0 \le i \le N, & p< j < N-p,\\
\tilde{\phi}_{p,i}(x) \tilde{\phi}_{p,j}(y), & i,j\leq p, &  i,j \geq N-p.
\end{array}
\right.
\end{equation}

Below, for notational convenience, we will not explicitly distinguish between the modified basis functions and the interior basis functions and simply drop the tilde notation. Also, we will use the notation $\mathbb{Q}_{p,N}$ to denote the one dimensional space spanned by the $(N+1)$ Galerkin difference basis functions associated with the nodal degrees of freedom. 

\subsubsection{Alternative Galerkin Difference Spaces}\label{DGDspace}

Lastly we note that instead of the locally (within element) continuous Galerkin basis, one can use discontinuous even-degree polynomials
constructed using cell-centered interpolation nodes as in \ci{GDDG,sipDGGD}. In this case we would need to introduce fluxes not only at the
boundaries of the macro-elements as we do here, but also at the boundaries of the cells within each macro-element. The construction of the
energy-based discontinuous Galerkin method would follow in the same way as presented below, and we could also construct the improved basis for
elements in regions where the wave speed is constant. We note that in this case it would be possible to make different choices for the fluxes for
the cell boundaries interior to each element and those between the macro-elements. However, here we will focus solely on the use of the
locally continuous basis functions.

\subsection{The Energy Based Discontinuous Galerkin Method for the Wave Equation}\label{GD_cost_subsection_1}
We consider a mesh that discretizes $\Omega$ into non-overlapping box shaped elements $\Omega^{k}$ with $\Omega = \bigcup_{k=1}^{m} \Omega^k$. On each element $\Omega^k = \otimes_{l \in \mathcal{S}}^d[L_{l}^k,R_{l}^k]$,  where $\mathcal{S} = \{x,y\}$ if $d=2$ and $\mathcal{S} = \{x,y,z\}$ if $d=3$. Let $(\mathbb{Q}_{p,N})^d$ be the space of functions spanned by the tensor product of the one dimensional Galerkin difference basis on an element. Then a test function $\varphi$ in $(\mathbb{Q}_{p,N})^d$ can be expressed (in three dimensions) as   
\begin{eqnarray*}
\varphi(x,y,z)= \prod_{\rho \in \mathcal{S}} \phi_{p,i_\rho}(\rho).
\end{eqnarray*}
Here we assume the same degree of approximation and the same number of degrees of freedom in each dimension, but remark that these can also be chosen independently.  

Now following \cite{DATH_UP} we define our discretization by the element-wise variational statement  
\begin{eqnarray} \label{DG1}
\int_{\Omega^k} c^2\nabla \varphi_u \cdot \left( \frac{\partial \nabla u}{\partial t}-\nabla v \right) d \Omega^k &=& \int_{\partial \Omega^k} c^2(\vec{n} \cdot \nabla \varphi_u) (v^\ast-v)
 \, d \Omega^k, \\
 \label{DG2}
\int_{ \Omega^k} \varphi_v \frac{\partial v}{\partial t}+c^2 \nabla \varphi_v \cdot\nabla u \, d \Omega^k &=&  \int_{\partial \Omega^k} c^2 \varphi_v \left( \vec{n} \cdot ( \nabla u)^\ast  \right) d \Omega^k,
\end{eqnarray}
for all $(\varphi_u,\varphi_v) \in (\mathbb{Q}_{p,N})^d \times (\mathbb{Q}_{p,N})^d$. As equation (\ref{DG1}) vanishes for constants we augment it by the independent equation 
\begin{equation}    
\int_{\Omega^k} \left( \frac{\partial u}{\partial t}-  v \right) d \Omega^k =0. \label{DG1_extra}
\end{equation}

\subsubsection{Numerical Fluxes}
Following the  notation in \cite{DATH_UP} we take $\beta,\tau\geq0$, $0\leq\alpha\leq 1$. Now, let the superscript $1$ represent data from inside an element and the superscript $2$ represent data from the outside of an element. Then we can write the numerical fluxes as
\begin{eqnarray*}
v^\ast &=& \alpha v^1+(1-\alpha)v^2-\beta(\nabla u^1\cdot {\bf n}^1 + \nabla u^2\cdot{\bf n}^2),\\
(\nabla u)^\ast &=& (1-\alpha)\nabla u^1+\alpha\nabla u^2-\tau(v{\bf n}^1+v{\bf n}^2).
\end{eqnarray*}
The above fluxes are energy-conserving when $\beta=\tau=0$, and upwind and dissipative when $\beta, \tau>0$. In the rest of the analysis, we focus on the three popular choices:
\begin{itemize}
\item[] {\it Central\ \ flux}: $\alpha = \frac{1}{2}$, $\beta = \tau = 0$,
\item[] {\it Alternating\ \ flux}: $\alpha = 1$, $\beta = \tau = 0$,
\item[] {\it Upwind\ \ flux}: $\alpha = \frac{1}{2}$, $\beta = \frac{\xi}{2}$, $\tau = \frac{1}{2\xi}$.
\end{itemize}
In the last flux $\xi$ is a flux splitting parameter with the same dimensional units as the speed of sound $c$.

We note that the possibility of choosing simple mesh-independent flux parameters is a feature of the energy-DG formulation.

\section{Efficient Formulation on Cartesian Grids}\label{GD_cost}
In this section we restrict our attention to the case of constant speed of sound within an un-mapped Cartesian element. Using the
tensor products of the compactly-supported basis functions described in Section~\ref{GD_basis} we will see that the complexity of
computing the time derivatives is unacceptable for $N$ large due to the structure of the lift matrix associated with (\ref{DG1}). 
However, if we use
a simple simultaneous diagonalization of the mass and stiffness matrix we can construct practical implementations of (\ref{DG1}), (\ref{DG2}) and (\ref{DG1_extra}) whose cost scales linearly with the total number of degrees of freedom.  

On each element $\Omega^k$, we approximate the solution by tensor product expansions (here for the case of three dimensions)
\begin{eqnarray*}
u(x,y,z,t) &=& \sum_{l_x=0}^{N} \sum_{l_y=0}^{N} \sum_{l_z=0}^{N} u_{l_x,l_y,l_z} \phi_{p,l_x}(x) \phi_{p,l_y}(y) \phi_{p,l_z}(z),\\
v(x,y,z,t) &=& \sum_{l_x=0}^{N} \sum_{l_y=0}^{N} \sum_{l_z=0}^{N} v_{l_x,l_y,l_z} \phi_{p,l_x}(x) \phi_{p,l_y}(y)\phi_{p,l_z}(z).
\end{eqnarray*}
On element $\Omega^k$, let $U^k$, $V^k$ be column vectors containing the nodal values of $u$ and $v$, that is $u_{l_x,l_y,l_z}$ and $v_{l_x,l_y,l_z}$, respectively. Then, in $d$-dimensions we may write the nodal based version of the method as the system of ordinary differential equations 
\begin{eqnarray}
\hat{S} \frac{d U^k}{dt} &=& \hat{S} V^k  \nonumber \\
&& + 
\sum_{j=1}^d 
(\alpha-1)\left[(D^{R,R}-D^{L,L})V^k
- D^{R,L}V^{k+1}
+ D^{L,R}V^{k-1}\right] \nonumber\\
&& - \sum_{j=1}^d  \beta \left[ (C^{R,R}-C^{L,L}) U^k
+ C^{R,L} U^{k+1}
 - C^{L,R}U^{k-1} \right], \label{DG_matrix1}
\end{eqnarray}
\begin{eqnarray}
M \frac{d V^k}{d t }  &=& - c^2 S U^k  \nonumber\\
&& + c^2\sum_{j=1}^{d}
(1-\alpha)\left( E^{R,R} -E^{L,L} \right)U^k
+ \alpha \left(E^{R,L}U^{k+1}
- E^{L,R}U^{k-1} \right) \nonumber\\
&&- c^2\sum_{j=1}^{d} \tau\left[ \left( B^{R,R} - B^{L,L} \right) V^k
+ B^{R,L}V^{k+1}
- B^{L,R}V^{k-1} \right].\label{DG_matrix2}
\end{eqnarray}
Here we abuse the notation in that for each coordinate direction in the sums we use the superscript $ k \pm 1$ denote the element ``left'' and ``right'' of element $k$ in the $j$th direction. The definitions of the mass matrix $M$, the stiffness matrix $\hat{S}$ and the lift matrices $B,C,D,E$ will be given below.

\subsection{Complexity with Galerkin Difference Basis}
Now, in order to compute the time derivatives $\frac{d U^k}{d t }$ and $\frac{d V^k}{d t}$ we must evaluate the matrix vector products on the right hand side and the action of the lift matrices on $U^k$ and $V^k$. As the matrices are sparse it is possible to do this at a cost that scales as $\sim f(p) N^\kappa$, with $f(p)$ being a low degree polynomial in $p$ and $\kappa = d$ for the volume terms and $\kappa = d-1$ for the surface terms. In other words the cost scales linearly with the number of degrees of freedom. 

Further, due to the tensor product structure of the mass matrix we have that the element mass matrix $M$ can be composed as a Kronecker product of the one dimensional matrix, which we denote $M_j$,   
\[
M = \otimes_{j=1}^{d}  M_j,
\]
with $M_{j,kl} = \int_{L_j}^{R_j} \phi_{j,k}\phi_{j,l}d_{x_j}$. Now, as the one dimensional mass matrices have bandwidth $p$ so will its $LU$-factors. Let $ L_j U_j \equiv M_j$. Then by the Hadamard product property we have  
\begin{equation*}
M = \otimes_{j=1}^d M_j = \otimes_{j=1}^d L_j U_j = (\otimes_{j=1}^d L_j)(\otimes_{j=1}^d U_j) \equiv LU.
\end{equation*}
Thus, as the cost of each substitution is $\mathcal{O}(p(N+1))$ the cost of solving $Mx = b$ is $\mathcal{O}((p(N+1))^d)$, which again is linear in the degrees of freedom. Unfortunately,  the stiffness matrix $\hat{S}$ is a sum of Kronecker products 
\begin{equation*}
\hat{S} = S_1\otimes M_2 \otimes M_3 + M_1\otimes S_2\otimes M_3 + M_1\otimes M_2\otimes S_3,  
\end{equation*}
with $S_{j,kl} = \int_{L_j}^{R_j}\frac{d\phi_{j,k}}{dx_j}\frac{d\phi_{j,l}}{dx_j}$, supplied by an extra equation from (\ref{DG1_extra}).
As is well-known, the cost of solving such a system directly is generally superlinear in the number of degrees-of-freedom, although recent
advances combining nested dissection ordering with low rank approximations can reduce this to near linear cost but with a significant prefactor
\ci{SchmitzYing3D}. The structure of (\ref{DG_matrix1}) can be exploited to rewrite it in the form:
\bd
\f {dU^k}{dt} = V^k + G (U^{k-1},V^{k-1},U^k,V^k,U^{k+1},V^{k+1})^T ,
\ed
where $G$ may be viewed as a lift matrix. As $G$ maps boundary data to volume data the cost of this operation will scale as
$\mathcal{O}((N+1)^{2d-1})$ per time step, with the larger cost of inverting $\hat{S}$ restricted to a precomputation. As $N$ may be much larger than $p$, this scaling implies that the resulting method would not be competitive with an implementation using a standard continuous Galerkin difference formulation of \cite{Gdiff} or the more recent method using a SIPG
formulation \ci{sipDGGD}. We must thus seek an improved method.

\subsection{Optimal Computational Complexity by Simultaneous Diagonalization}\label{GD_new_basis}
The above mentioned complexity for evolving $U$ is not competitive for practical computations unless we limit the element sizes, e.g.
with $N = \mathcal{O} (p^{(d-1)/d})$. In this section we follow  \cite{lynch1964direct} and show that it is possible to make a simple (computational) change of basis that results in a method with linear complexity. Precisely the new basis is found by solving the generalized eigenvalue of problem for each of the one dimensional matrices, $S_j$, $j = 1,\ldots,d$. That is, the new basis vectors are solutions to,
\begin{eqnarray*}
S_j\psi_{j,k_j} = \lambda_{j,k_j} M_j\psi_{j,k_j}, \ \ k_j=0,\ldots,N_j, \ \ j=1,\ldots,d.
\end{eqnarray*} 
We normalize the eigenvectors according to
\begin{eqnarray*}
\psi_{j,k_j} \leftarrow \frac{\psi_{j,k_j}}{((\psi_{j,k_j})^TM_j\psi_{j,k_j})^{1/2}}.
\end{eqnarray*}

Let $\Psi_j $ be the matrix containing the new one dimensional basis 
\[
\Psi_j = \left(
\begin{array}{cccc}
\psi_{j,0} & \psi_{j,1} & \cdots & \psi_{j,N_j}
\end{array} 
\right) .
\]
Then the $d$-dimensional basis is 
\begin{equation*}
\Psi = \otimes_{j=1}^d \Psi_j.
\end{equation*}
Now, we define $\bar{U}$ and $\bar{V}$ by 
\[
U = \Psi \bar{U}, \ \ V = \Psi \bar{V},
\]
then equation (\ref{DG1}) and (\ref{DG2}) become 
\begin{eqnarray}
&&\Psi^TS\Psi \left( \frac{d \bar{U}}{dt}-\bar{V} \right) = f_u, \label{DG1_sparse} \\
&&\Psi^TM\Psi \frac{d \bar{V}}{dt} +c^2\Psi^TS\Psi\bar{U}= f_v. \label{DG2_sparse}
\end{eqnarray}
where,
\begin{eqnarray*}
f_u &=& \sum_{j=1}^d 
(\alpha-1) \Psi^T \left[(D^{R,R}-D^{L,L})\Psi \bar{V}^k
- D^{R,L}\Psi \bar{V}^{k+1}
+ D^{L,R}\Psi \bar{V}^{k-1}\right] \\
&& - \sum_{j=1}^d  \beta \Psi^T \left[ (C^{R,R}-C^{L,L}) \Psi \bar{U}^k
+ C^{R,L} \Psi \bar{U}^{k+1}
 - C^{L,R}\Psi \bar{U}^{k-1} \right],\\
f_v &=&c^2 \sum_{j=1}^{d}
(1-\alpha) \Psi^T \left( E^{R,R} -E^{L,L} \right)\Psi \bar{U}^k
+ \alpha \Psi^T  \left(E^{R,L}\Psi \bar{U}^{k+1}
- E^{L,R}\Psi \bar{U}^{k-1} \right) \\
&&- c^2\sum_{j=1}^{d} \tau \Psi^T\left[ \left( B^{R,R} - B^{L,L} \right) \Psi \bar{V}^k
+ B^{R,L}\Psi \bar{V}^{k+1}
- B^{L,R}\Psi \bar{V}^{k-1} \right].
\end{eqnarray*}

In the new basis we have that the mass matrix diagonalizes 
\begin{equation*}
\Psi^T M \Psi =
\otimes_{j=1}^d (\Psi_j^T M_j \Psi_j)  = I,
\end{equation*}
as does the differentiation matrix 
\begin{multline*}
\Psi^T S \Psi = \sum_{j=1}^d 
(\Psi_1^T M_1 \Psi_1) 
\otimes \cdots \otimes 
(\Psi_{j-1}^TM_{j-1}\Psi_{j-1}) 
\otimes (\Psi_j^TS_j\Psi_j)\otimes(\Psi_{j+1}^TM_{j+1}\Psi_{j+1})\otimes 
\cdots \\
\cdots
\otimes (\Psi_d^T M_d \Psi_d) =  \sum_j \Lambda_j,
\end{multline*}
where,
\begin{equation*}
    \Lambda_j = I
    \otimes
    \cdots
    \otimes I
    \otimes
    \text{diag}(\lambda_{j,0},\ldots,\lambda_{j,N_j})
    \otimes I
    \otimes
    \cdots
    \otimes I.
\end{equation*}
In the above equations we use $I$ to denote the identity matrix of size inferred by the context.

Note that as one of the eigenvalues of $S$ is zero one of the equations in (\ref{DG1_sparse}) vanishes. Suppose we have ordered the unknowns so that this corresponds to the first entry in $\bar{U}$, then we simply enforce the additional independent equation 
\[
\frac{d \bar{U}_1}{dt} = \bar{V}_1.
\]   

We thus conclude that the cost of all the volume terms scales linearly with the number of degrees of freedom. We now turn to the evaluation of the surface terms in the new basis. 

Consider first the surface terms in the Galerkin difference basis. In a single dimension the elements in the four different surface terms are of the form
\begin{equation*}
\tilde{B}_{kl}^{X,Y} = \phi_{k}(X)\phi_{l}(Y), \ \ 
\tilde{D}_{kl}^{X,Y} = \frac{d\phi_{k}}{dx}(X)\phi_{l}(Y), \ \ 
\tilde{E}_{kl}^{X,Y} = \phi_{k}(X) \frac{d\phi_{l}}{dx}(Y), \ \ 
\tilde{C}_{kl}^{X,Y} = \frac{d\phi_{k}}{dx}(X)\frac{d\phi_{l}}{dx}(Y),
\end{equation*}
where $\{X,Y\}\in\{\{L,L\}, \{L,R\}, \{R,L\}, \{R,R\}\}$. 

Now due to the local support properties of the Galerkin difference basis the number of nonzero elements in the above matrices are $1$ for $\tilde{B}^{X,Y}$, $(p+1)$ for $\tilde{D}^{X,Y}$ and $\tilde{E}^{X,Y}$, and $ (p+1)^2$ for $\tilde{C}^{X,Y}$. 

The $d$-dimensional version of the surface matrices can again be constructed by Kronecker products. For example we have that
\begin{eqnarray*}
D_j^{X,Y} = M_1\otimes\cdots\otimes M_{j-1}\otimes\tilde{D}_j^{X,Y}\otimes M_{j+1}\otimes\cdots\otimes M_{d}.
\end{eqnarray*}

Applying the change of basis we have that 
\begin{multline*}
\Psi^T D_j^{X,Y} \Psi = (\Psi_1^T M_1 \Psi_1)\otimes\cdots\otimes(\Psi_{j-1}^T M_{j-1}\Psi_{j-1})\otimes(\Psi_j^T \tilde{D}_j^{X,Y}\Psi_j)\otimes(\Psi_{j+1}^TM_{j+1}\Psi_{j+1})\otimes \cdots \\ 
\cdots \otimes(\Psi_d^T M_d \Psi_d)
= I\otimes\cdots\otimes I\otimes(\Psi_j \tilde{D}_j^{X,Y} \Psi_j)\otimes I\otimes\cdots\otimes I.
\end{multline*}
Thus applying $\Psi^T D_j^{X,Y} \Psi$ to the $(N+1)^d$ dimensional vector $\bar{V}$ can be done at a cost that scales with $ (p+1) (N+1)^d$. Similarly, the cost of applying   $\Psi^T B_j^{X,Y} \Psi$, $\Psi^T C_j^{X,Y} \Psi$, and $\Psi^T E_j^{X,Y} \Psi$ can be done at a cost of $ (N+1)^d$, $ (p+1)^2 (N+1)^d$, and $ (p+1) (N+1)^d$, respectively. 

\graphicspath{{graphs_computational_cost/}}
\begin{figure}[htb]
\begin{center}
\includegraphics[width=0.45\textwidth]{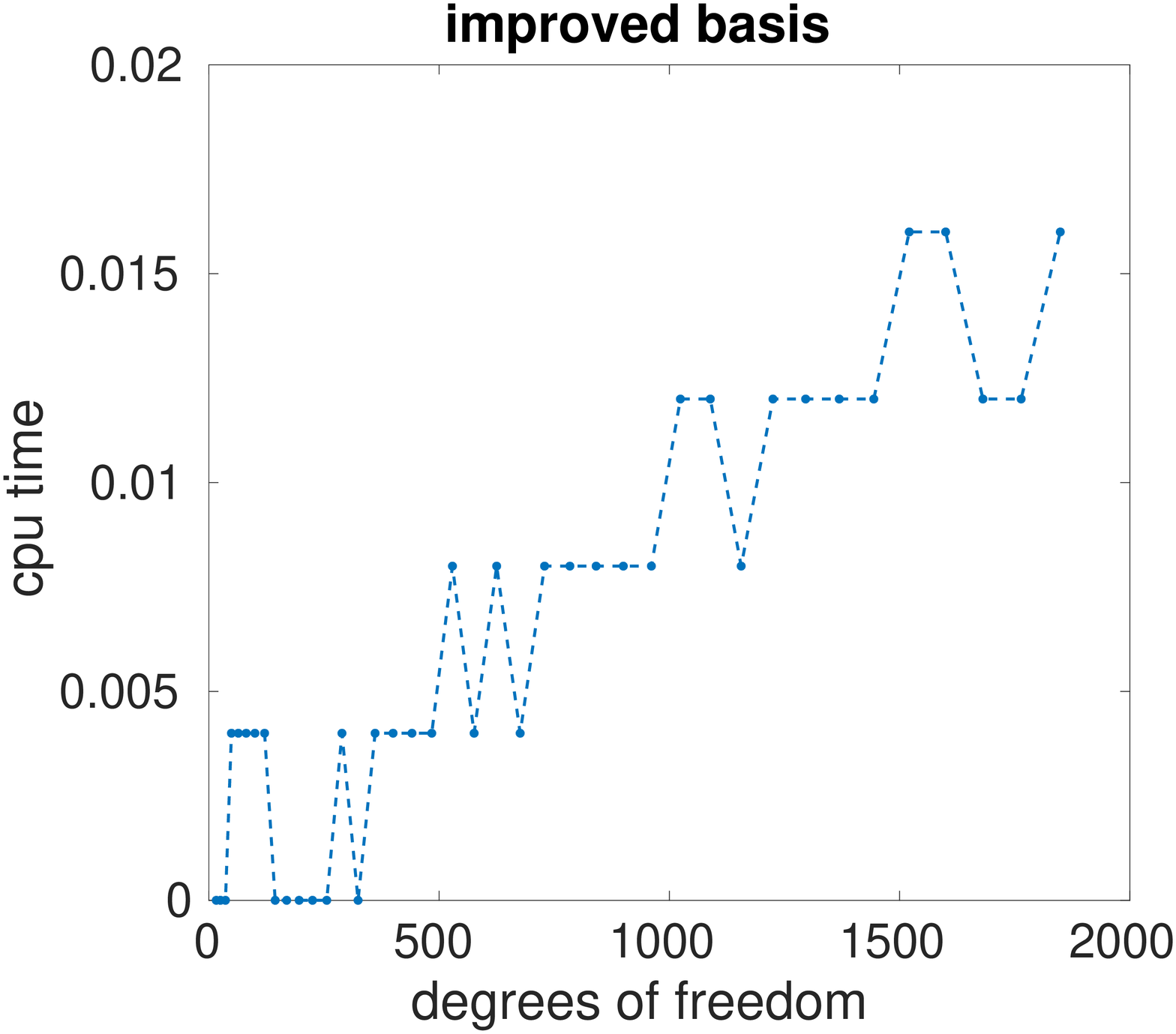}
\includegraphics[width=0.45\textwidth]{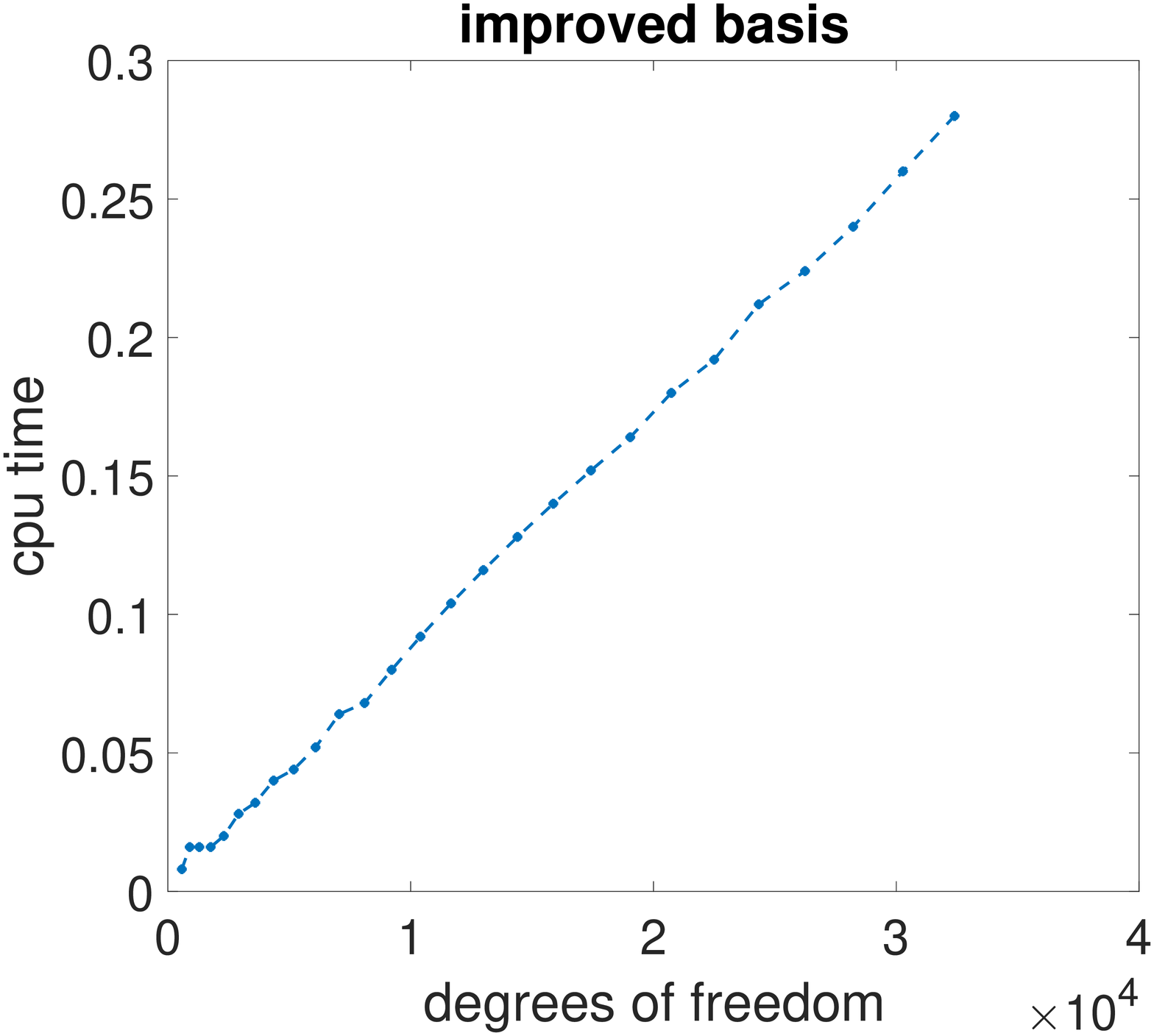}
\caption{Plots of CPU time (in seconds) in the Galerkin difference basis obtained by simultaneous diagonalization as a function of degrees of freedom (DOF) in two space dimensions. The graph on the left is for one DG element $(n = 1)$ and on the right for 36 DG elements $(n=6)$. The speed of sound is $c = 1$ and the splitting parameter for the upwind flux is $\xi = c$.\label{fig:timing_2d}} 
\end{center}
\end{figure}

\subsection{Numerical Verification of the Computational Complexity}
We now present timing results that illustrate above analysis. We consider a two dimensional problem in the domain $[0,1]\times[0,1]$ and use the upwind flux (the other fluxes give similar timing results). For this problem we choose forcing and boundary conditions so that the solution is
\begin{equation*}
u(x,t) = \sin(16\pi t)\sin(16\pi x) + \cos(16\pi t)\cos(16\pi y).
\end{equation*}
We present results for both one DG element ($n=1$) and 36 DG elements ($n=6$). The degrees of freedom are $N+1$ for both the $x$ and the $y$ direction, i.e, the total degrees of freedom is $ n^2 (N+1)^2$. Specifically, we choose $N = 3,4,5\cdots,42$ when $n = 1$ and $N = 3,4,5\cdots,29$ when $n = 6$ for the improved basis. We also choose $p = 3$ which results in a fourth order accurate method.

To time the code we use the built-in function \verb|CPU_TIME()| in \verb|FORTRAN| to record the elapsed CPU time which is used to evolve the solution using the classic fourth order acurate Runge-Kutta method for $10$ timesteps steps. 

In Figure \ref{fig:timing_2d}, we observe that the CPU time is proportional to $n^2(N+1)^2$ for the Galerkin difference basis obtained by the simultaneous diagonalization, which is orders of magnitude smaller than what would be required for the standard Galerkin basis.

\section{Dispersion Analysis}\label{dispersion}
To investigate how well the scheme proposed in Section \ref{GD_cost_subsection_1} preserves the wave propagation properties, we use the standard Bloch wave analysis as in \cite{ainsworth2004dispersive,AinsworthMonkMuniz,hu1999analysis,moura2015linear,hesthaven2007nodal}.  Here, we consider
\begin{equation}\label{wave_1d_problem}
\frac{\partial^2 u}{\partial t^2} = c^2\frac{\partial^2 u}{\partial x^2}, \ \ a<x<b,
\end{equation}
with initial condition $u(x,0) = e^{ikx}$ and periodic boundary conditions $u(a,t) = u(b,t)$, $\f {\pa u}{\pa x}(a,t)= \f {\pa u}{\pa x}(b,t)$.
We then seek spatially periodic solutions of the form
\begin{equation*}\label{solution_dispersion}
u(x,t) = e^{i(\kappa x-wt)},
\end{equation*}
from which the exact dispersion relation $\omega = \pm c\kappa$ for (\ref{wave_1d_problem}) can be found.

Next, partition the computational domain into non-overlapping uniform DG elements $I^k = [x^k,x^{k+1}]$, $k = 0,\cdots,n-1$ with $H = x^{k+1}-x^k = (b-a)/n$. For each DG element $I^k$, there are $N+1$ equidistant nodal degrees of freedom with spacing $h = H/N$. The semi-discretization in $I^k$ becomes then follows from (\ref{DG_matrix1}) and (\ref{DG_matrix2}) with $d=1$. Let the vectors $U^k = (U_0^k,U_1^k,\cdots,U_{N}^k)$ and $V^k = (V_0^k,V_1^k,\cdots,V_{N}^k)$ hold the nodal approximations of $u(x,t)$ and $v(x,t)$ in $I^k$, respectively. 

We seek solutions in terms of Bloch waves
\begin{equation}\label{dispersion_transfer}
U^k_l = \hat{U}^k_le^{i(\kappa (x^k+lh)-\omega t)},\quad V^k_l = \hat{V}_l^ke^{i(\kappa (x^k+lh)-\omega t)},\ \ l=0,\cdots,N,
\end{equation}
and thus assume periodicity of the solution. That is 
\begin{equation}\label{periodic_dispersion}
W^{k+1}_0 = e^{i\kappa H}W^k_0, \quad W^{k-1}_N = e^{-i\kappa H}W^k_N,
\end{equation}
where $W$ represents either $U$ or $V$. To condense the notation, we have omitted the superscript $k$ for the rest of this section. Let $Z = (U,V)^T$ and combine with (\ref{DG_matrix1})-(\ref{DG_matrix2}), (\ref{dispersion_transfer})-(\ref{periodic_dispersion}) to obtain the following eigenvalue problem
\begin{equation*}
\hat{A}Z = -i\Omega Z,\ \ \hat{A} = \hat{A}_1\hat{A}_2,
\end{equation*}
where
\begin{equation*}
\Omega = \frac{\omega H}{c},\quad \hat{A}_1 = \frac{H}{c} 
\begin{pmatrix}
\hat{S}^{-1}&~\\
~&M^{-1}
\end{pmatrix},\ \
\hat{A}_2 = 
\begin{pmatrix}
A_2^{11}&A_2^{12}\\
A_2^{21}&A_2^{22}
\end{pmatrix},
\end{equation*}
with
\begin{eqnarray*}
	A_2^{11} &=& \tau (C^{L,L} - C^{R,R})-\tau (e^{iK}C^{R,L}-e^{-iK}C^{L,R}),\\
	A_2^{12} &=& \hat{S}+ (-\theta)(D^{R,R}-D^{L,L})+\theta(e^{iK}D^{R,L}-e^{-iK}D^{L,R}),\\
	A_2^{21} &=&-c^2S+ \theta c^2(E^{R,R}-E^{L,L})+(1-\theta) c^2(e^{iK}E^{R,L}-e^{-iK}E^{L,R}),\\
	A_2^{22}&=&-\beta c^2(B^{R,R}-B^{L,L})-\beta c^2(e^{iK}B^{R,L}-e^{-iK}B^{L,R}),
\end{eqnarray*}
and $K = \kappa H$. Note that the values of $\Omega = \Omega_r + i\Omega_i$ are, in general, complex valued. Here $\Omega_i$ is non-positive and represents the numerical damping of the corresponding scheme, and the real part $\Omega_r$ is the approximation to  $\frac{\omega H}{c}$.

For the numerical simulations in this section, the computational domain is chosen to be $x\in[0,1]$. The order of the approximation space is set to be $p = 3$. The number of degrees of freedom in each DG element is $N+1 = 10$. 
\graphicspath{{graphs_dispersion/}}
\begin{figure}[htb]
	\begin{center}
		\includegraphics[width=0.49\textwidth]{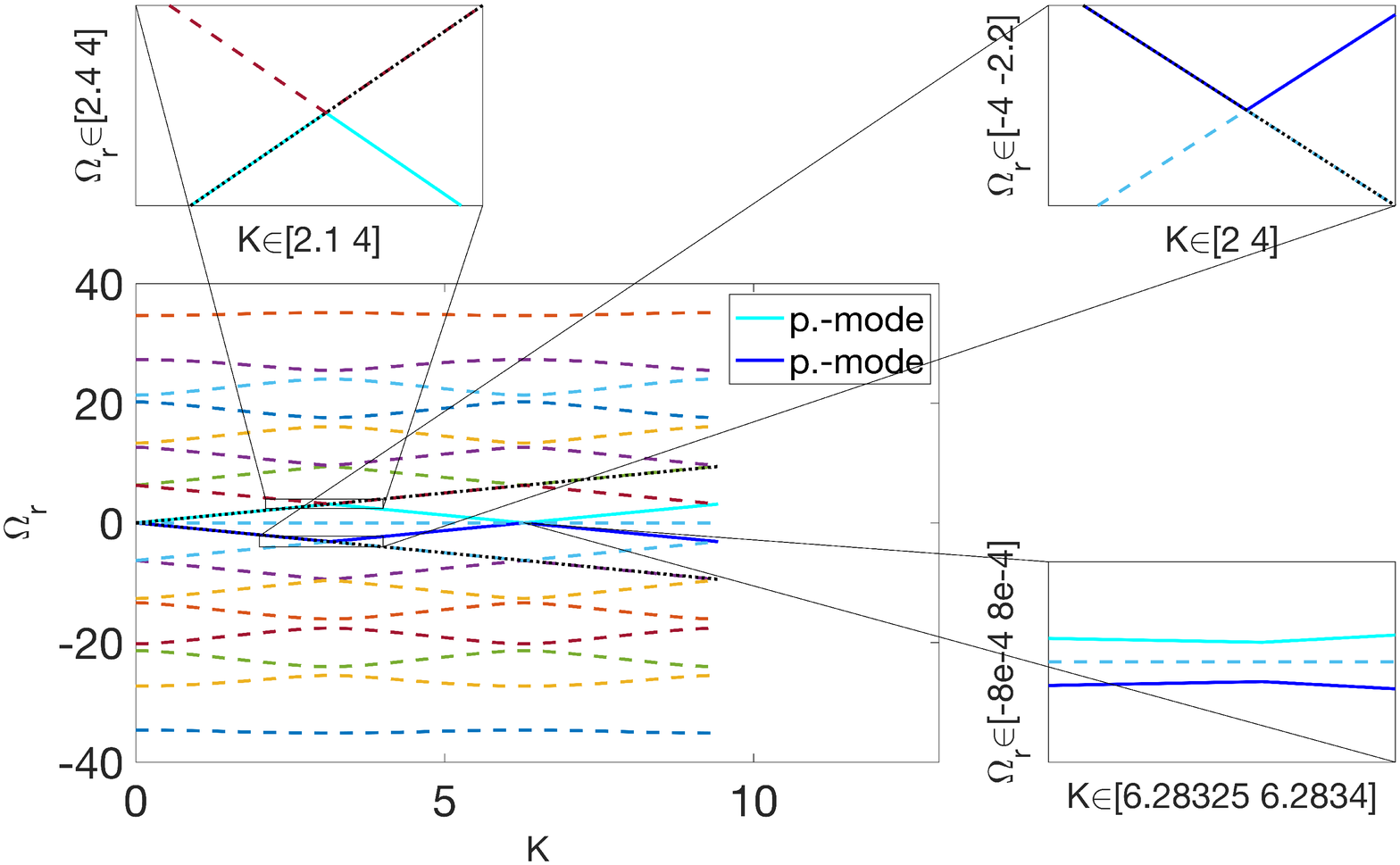}
		\includegraphics[width=0.49\textwidth]{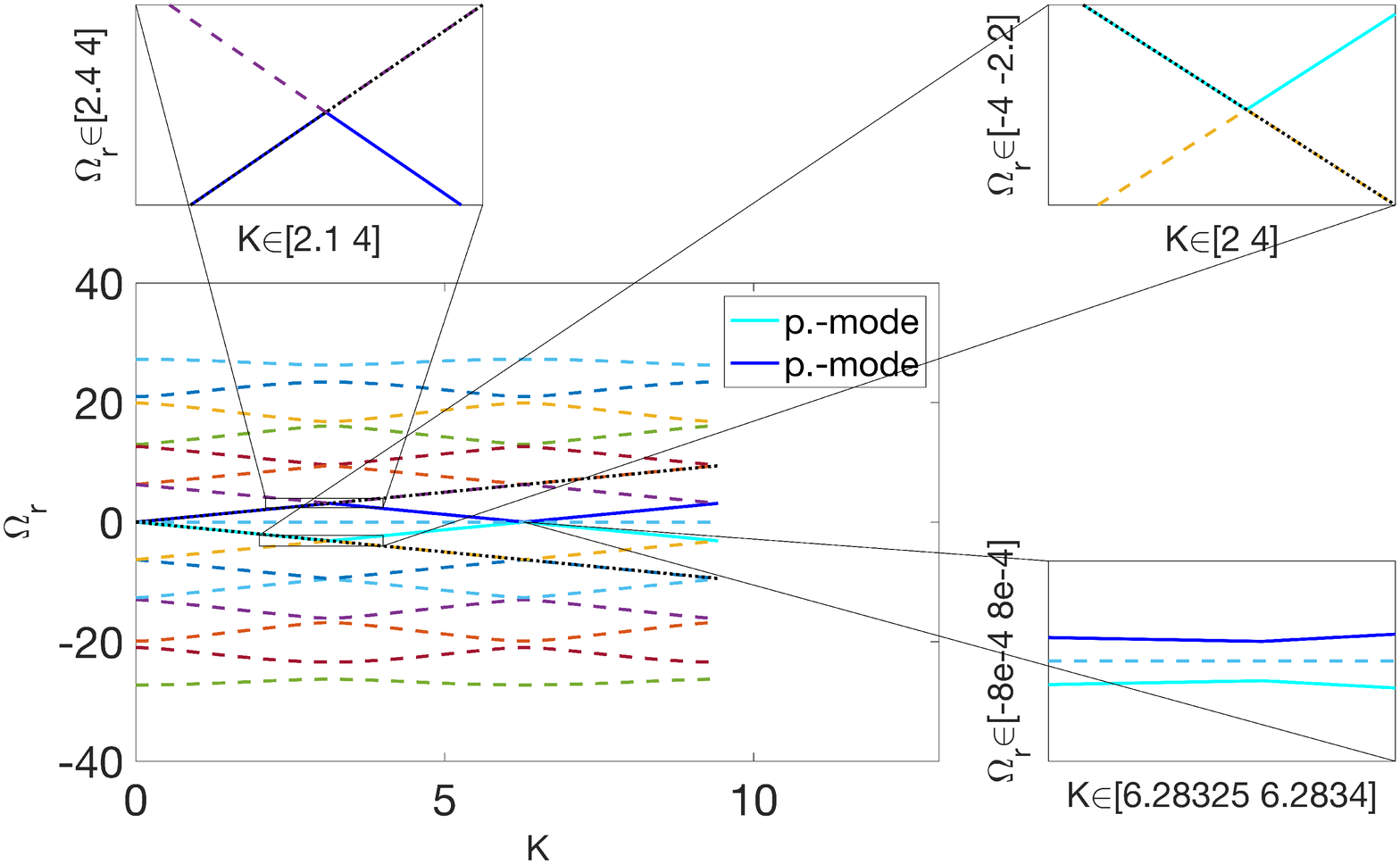}
	\end{center}
	\caption{In this simulation, we choose $n = 2$ DG elements and $N = 9$ GD cells in each DG element. On the left we show the numerical dispersion relation for the central flux. On the right we show the numerical dispersion relation for the alternating flux. Black dot lines are the exact dispersion relation, p.-mode represents physical modes, dashed lines are for spurious modes.\label{dispersion_dissipation_c_2ele}}
\end{figure}
\begin{figure}[htb]
	\begin{center}
		\includegraphics[width=0.49\textwidth]{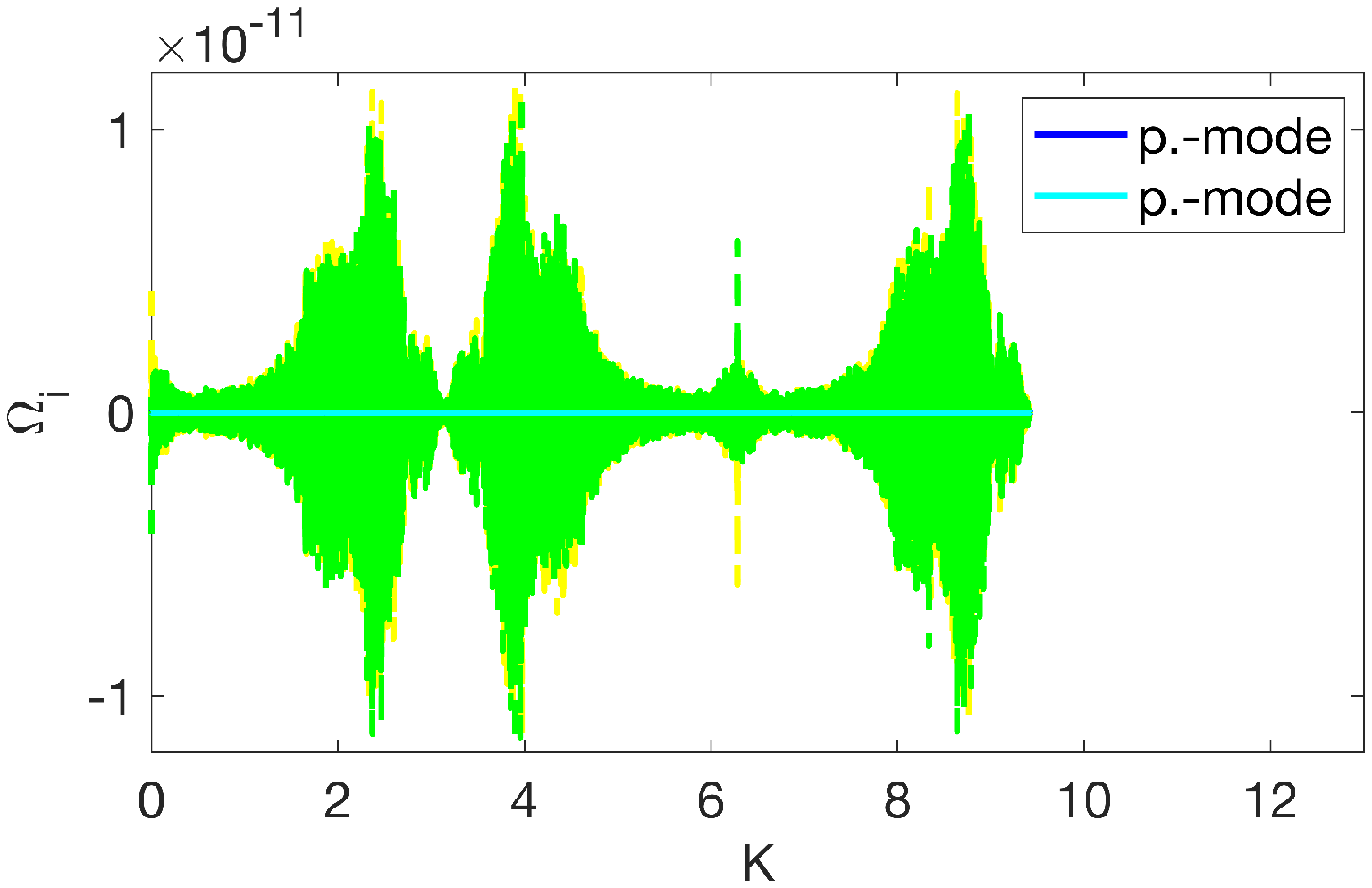}
		\includegraphics[width=0.49\textwidth]{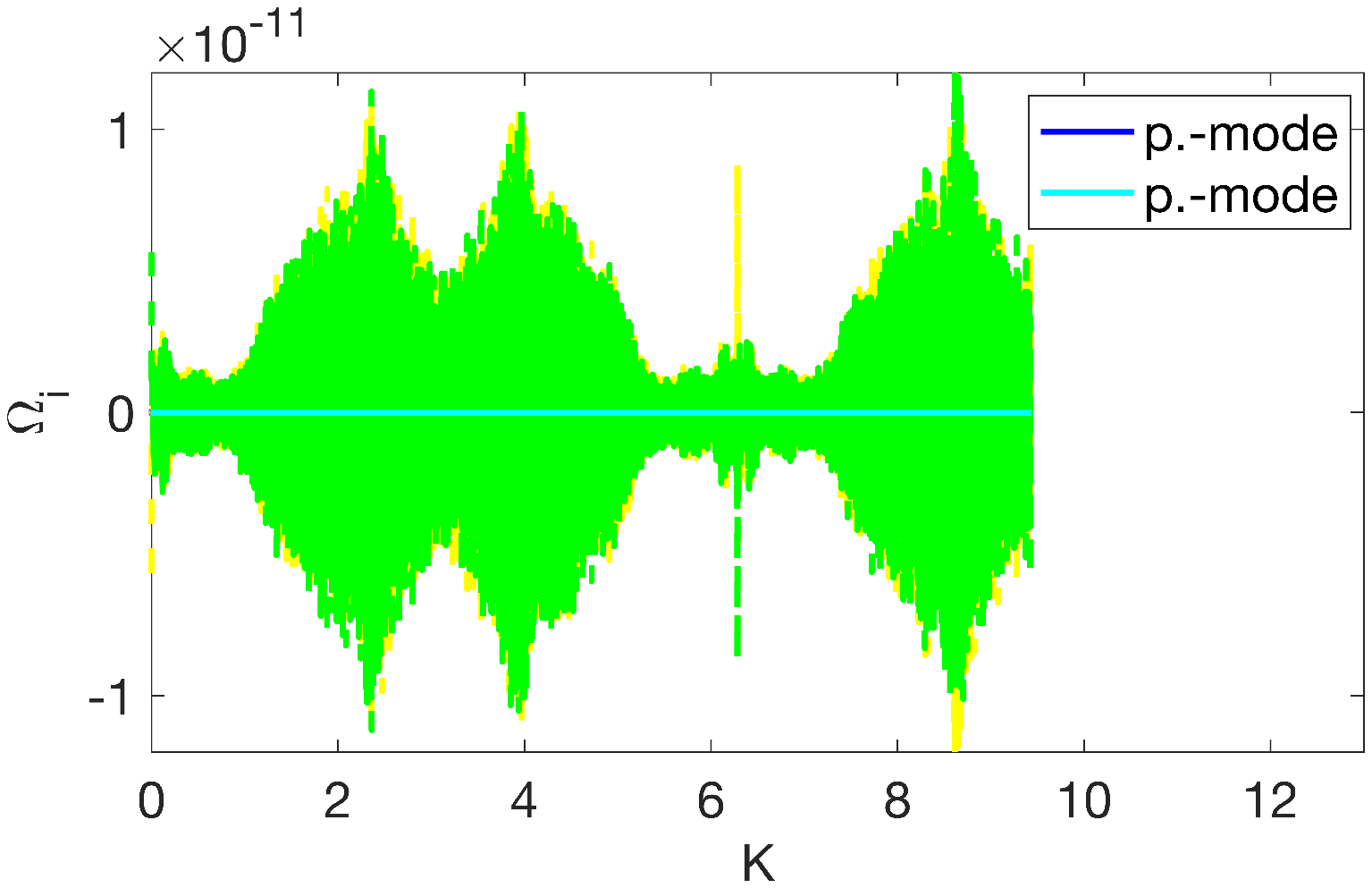}
	\end{center}
	\caption{In this simulation, we choose $n = 2$ DG elements and $N = 9$ GD cells in each DG element. On the left we show the numerical dissipation for the central flux for all modes. On the right we show the numerical dissipation for the alternating flux for all modes. 'p.-mode' represents the physical, dashed lines are for spurious modes.\label{dispersion_dissipation_a_2ele}}
\end{figure}

\begin{figure}[htb]
	\begin{center}
		\includegraphics[width=0.49\textwidth]{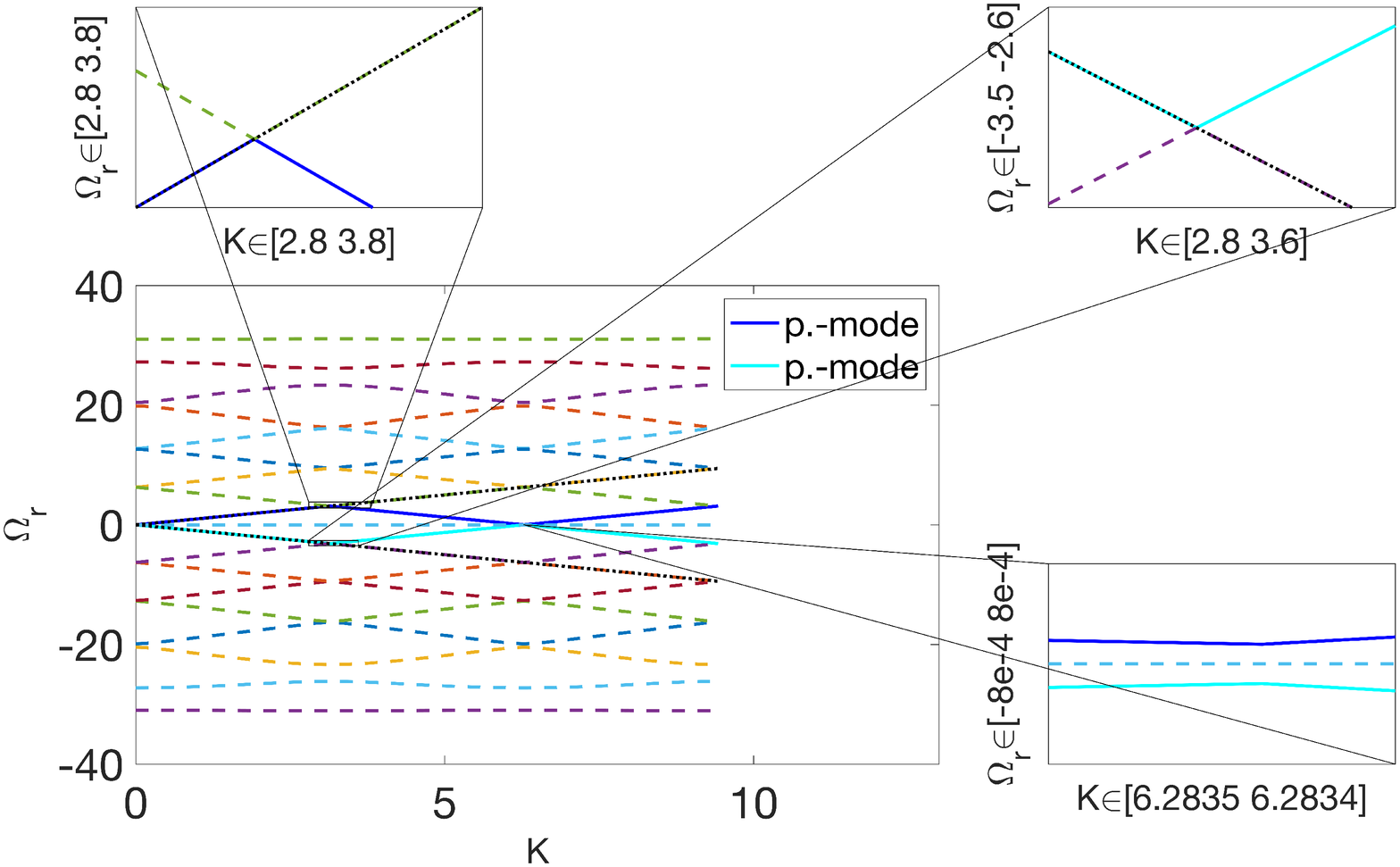}
		\includegraphics[width=0.49\textwidth]{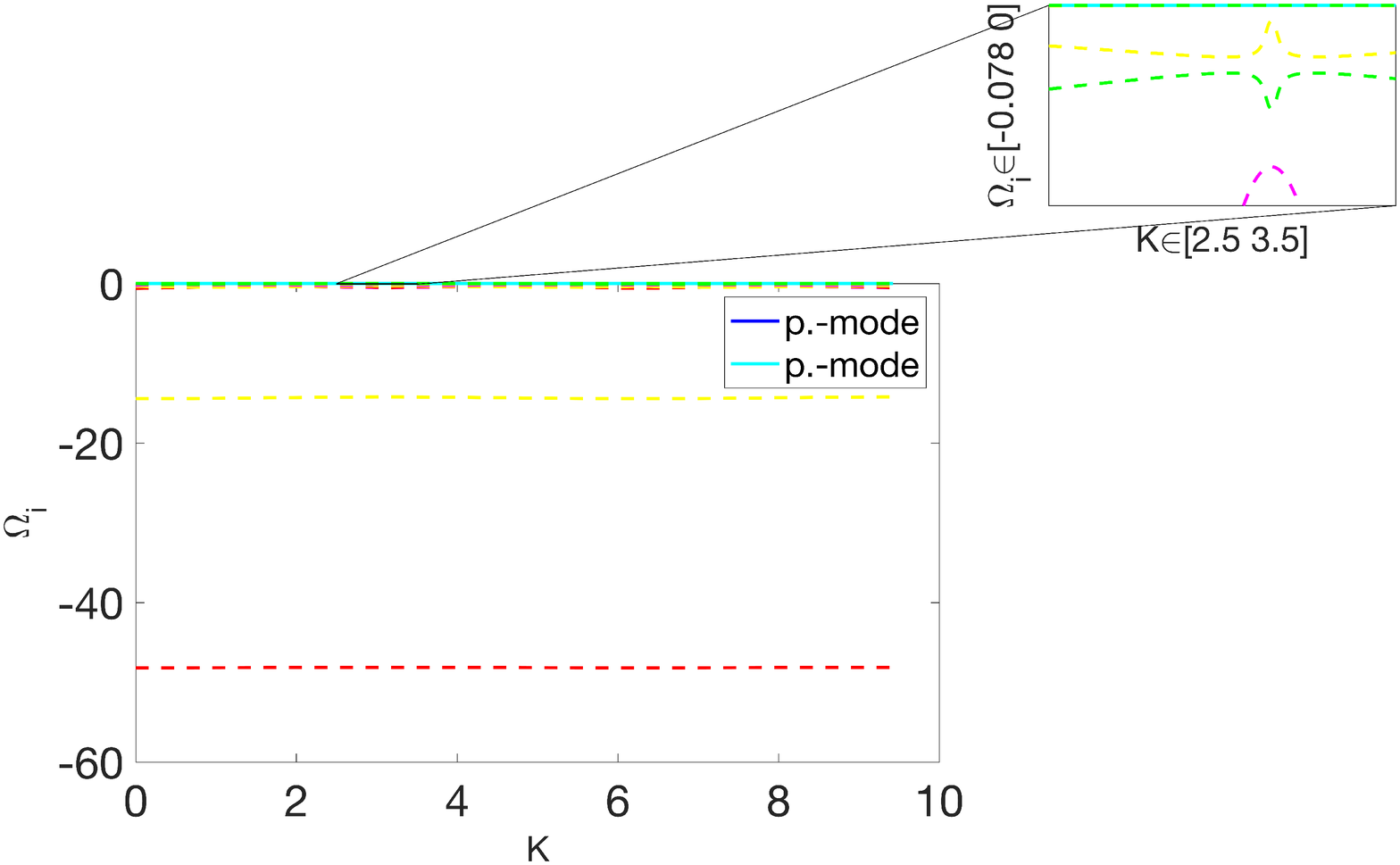}
	\end{center}
	\caption{In this simulation, we choose $n = 2$ DG elements and $N = 9$ GD cells in each DG element. On the left, we show the numerical dispersion relation for the upwind flux. Black dot lines are the exact dispersion relation,  p.-mode represents the physical and dash lines are for spurious modes. On the right, we illustrate the dissipation associated with the twenty modes.\label{dispersion_dissipation_up_2ele}}
\end{figure}

Figure \ref{dispersion_dissipation_up_2ele} presents the dispersion relation of the upwind flux. When $K$ is small, the numerical phase velocity also reflects the physical wave speed. Comparing the results in Figure \ref{dispersion_dissipation_up_2ele} and Figure \ref{dispersion_dissipation_c_2ele}--\ref{dispersion_dissipation_a_2ele}, we find that both conservative (central flux and alternating flux) and dissipative schemes (upwind flux) recover the physical mode for small $K$. The conservative schemes admit more complicated phenomena: the spurious modes do not damp for small values of $K$; for the dissipative scheme, however, the unphysical modes are strongly damped.

\begin{figure}[htb]
	\begin{center}
		\includegraphics[width=0.49\textwidth]{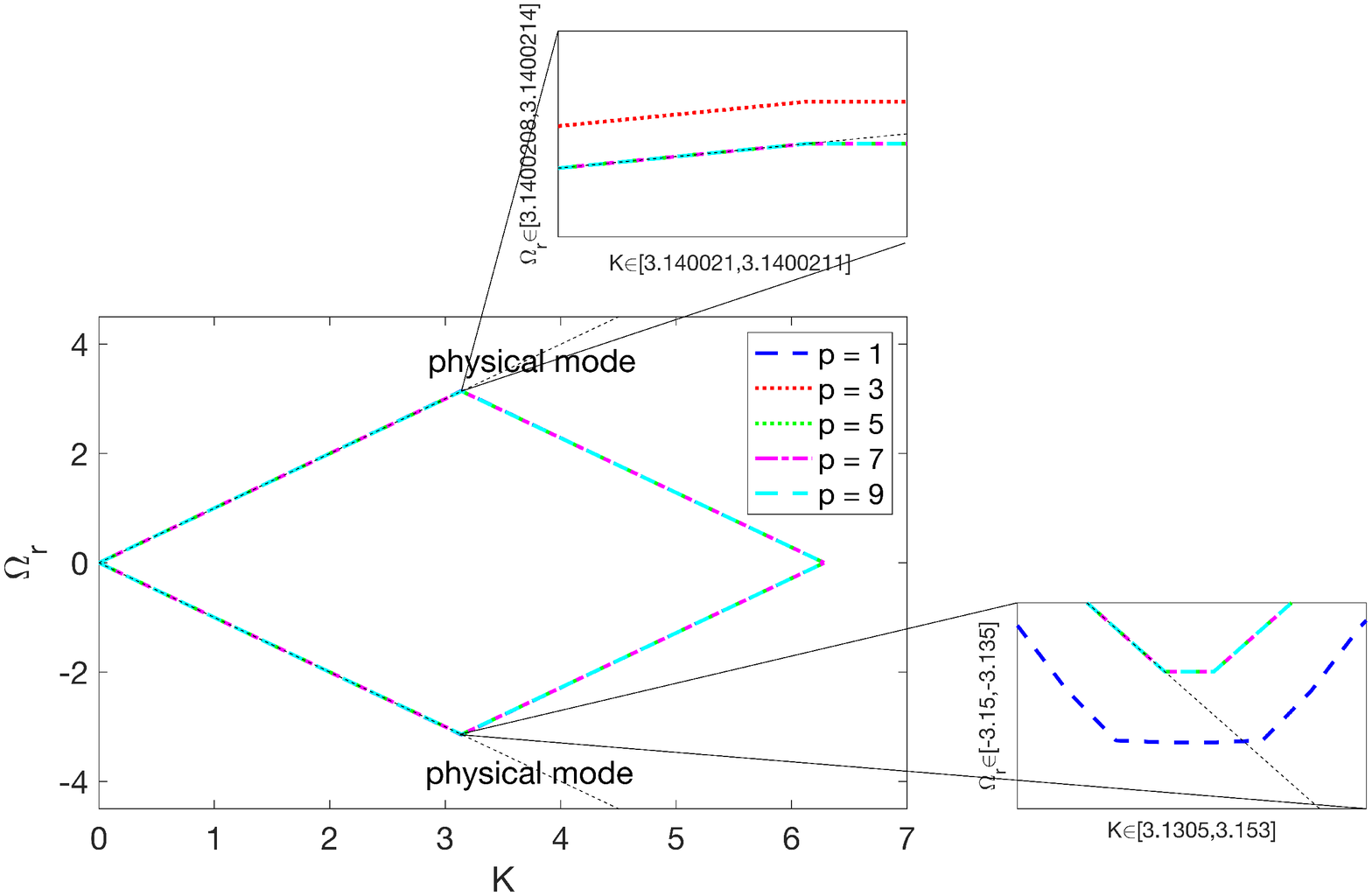}
		\includegraphics[width=0.49\textwidth]{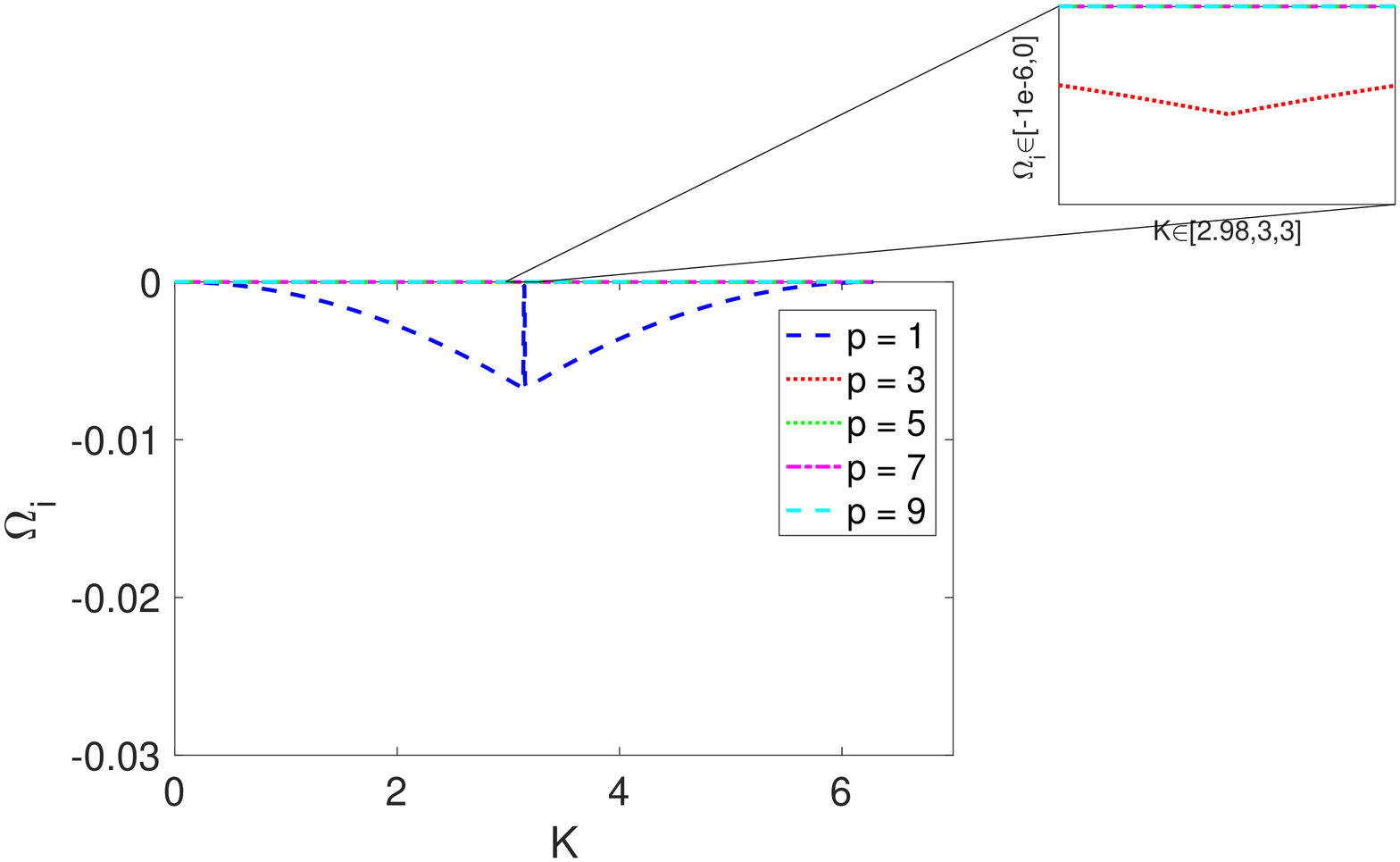}
	\end{center}
	\caption{On the left, we show the numerical dispersion relation for the upwind flux. The black dashed line is the exact dispersion relation. We present the numerical dispersion relation for the physical mode at different orders $p = (1,3,5,7,9)$ with $N = 19$. On the right, we illustrate the dissipation associated with the different orders $p = (1,3,5,7,9)$ with $N = 19$.\label{dispersion_dissipation_up_dq}}
\end{figure}
In Figure \ref{dispersion_dissipation_up_dq}, we show the dispersion relation of physical modes of the dissipative scheme (upwind flux) for a range of orders of approximation $p = (1,3,5,7,9)$ and $N = 19$ for all different approximation degrees $p$. We see that the numerical phase velocity is very close to the physical wave speed when $K$ is small and improves for a broader range of $K$ as the order of the approximation increases.

\begin{figure}[htb]
\begin{center}
\includegraphics[width=0.49\textwidth]{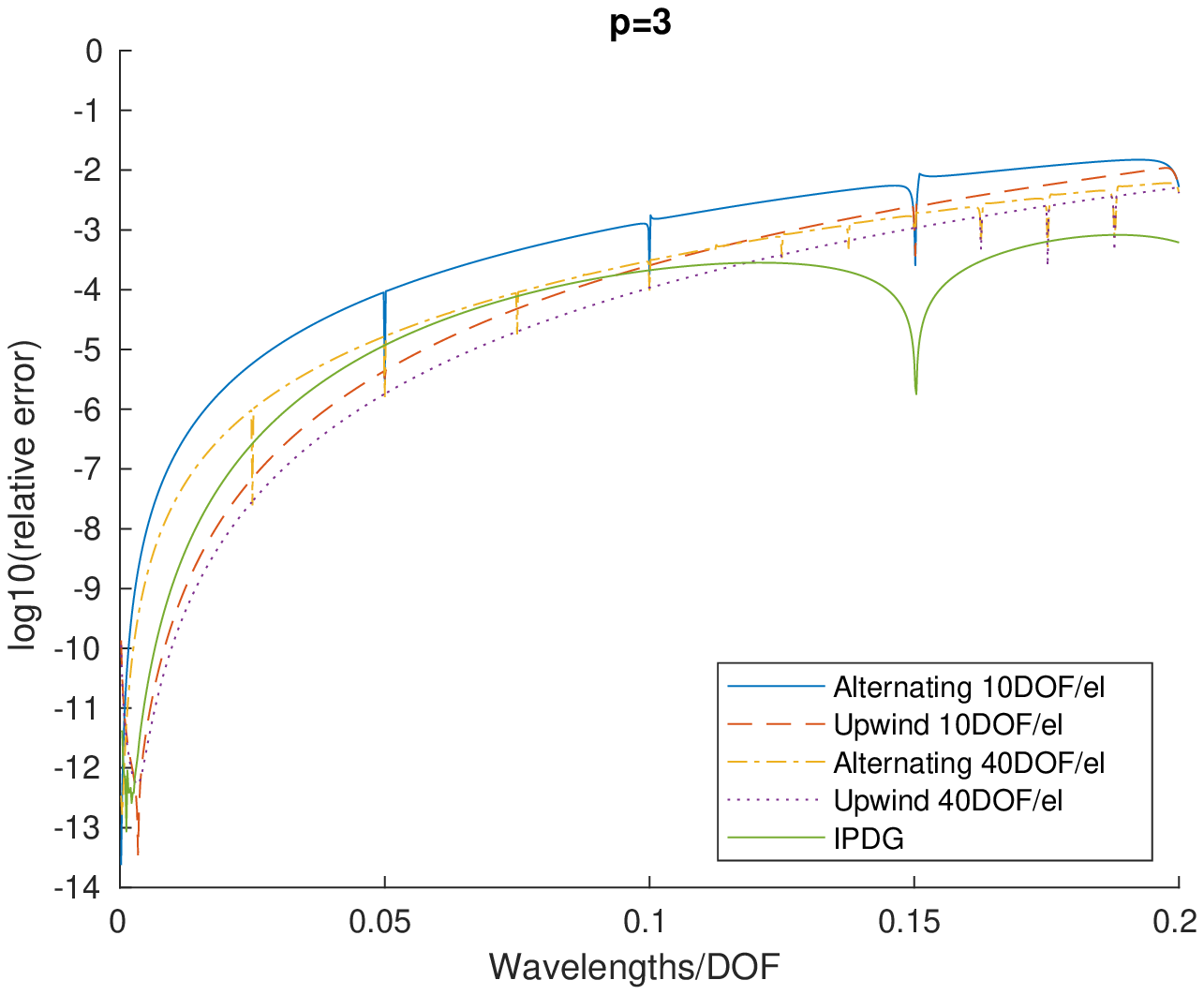}
\includegraphics[width=0.49\textwidth]{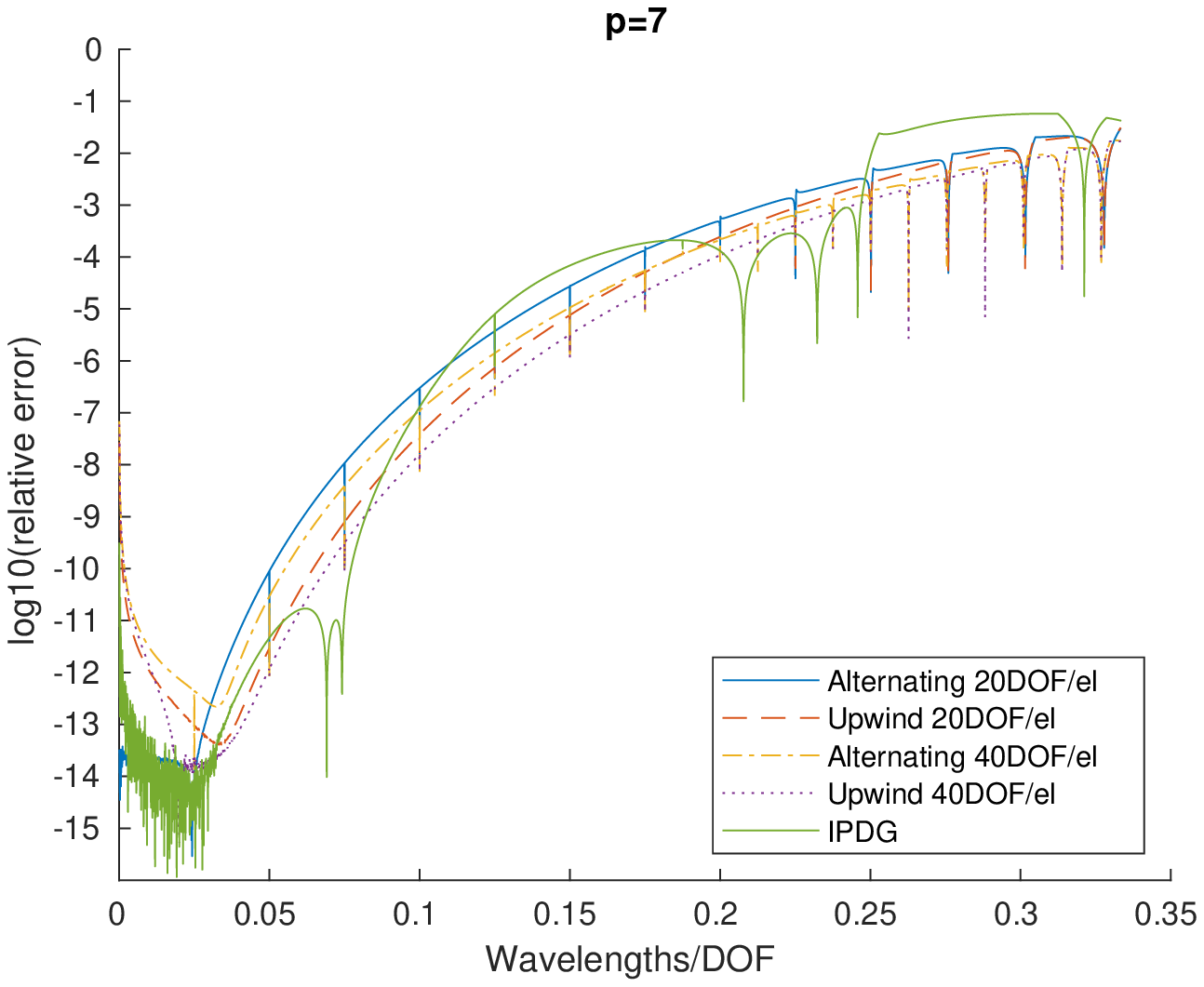}
\end{center}
\caption{On the left we compare the dispersion errors with $p=3$
for the alternating and upwind fluxes with $10$ and $40$ points per element as well as for IPDG scheme. On the right we show
the same for $p=7$ where the alternating and upwind fluxes are considered for $20$ and $40$ points per element.\label{discomp}}
\end{figure}

Lastly we consider direct comparisons of the proposed method with both alternating and upwind flux choices to the
interior penalty discontinuous Galerkin method (IPDG). We also consider the effect of increasing the number of nodes
within the element. The dispersion relation for the IPDG method is discussed extensively in \ci{AinsworthMonkMuniz}, and we
use the Bloch wave formulation and penalty parameter suggested there. On the left in Figure \ref{discomp} we take $p=3$ and
apply the method proposed here with $9$ and $39$ cells per element along with IPDG. In all cases we use a single
element and compute the dispersion
relation for a range of wave numbers with the largest wave number corresponding to $5$ degrees-of-freedom per wavelength. Note
that this implies a significant disparity in the wave numbers considered. Here we see that the dispersion error for the upwind
method is significantly smaller than for the conservative alternating flux scheme, but keep in mind that the graph does not
take the dissipation error of the upwind discretization into account. Both our conservative method and IPDG have oscillatory
dispersion relations, but perhaps because of the larger wave number range involved we see more oscillations for the proposed
method. Increasing the number of nodes per element improves the results. Here IPDG is the best performer at the coarsest
discretization level, though we do recall that the proposed scheme admits larger time steps than IPDG.

On the right we display the analogous plot for $p=7$ and largest wave number corresponding to $3$ degrees-of-freedom per
wavelength. Now the results for the proposed methods are shown with $19$ and $39$ cells per element. The comparisons between
them follow the same pattern as in the case $p=3$. However in this case the dispersion errors are larger for IPDG at the
coarse resolutions.

\section{Spectral Radius}\label{spectral_radius}
In this section, we study the spectral radius of the semi-discretization of our scheme. Consider the problem in one dimension
\[
\frac{\partial^2 u}{\partial t^2} = c^2\frac{\partial^2 u}{\partial x^2},\quad x\in[a,b],
\]
with either periodic boundary conditions or a homogeneous Dirichlet boundary condition at the left boundary and a homogeneous Neumann boundary condition at the right boundary. We emphasize that an expected advantage of Galerkin difference methods compared with standard DG schemes is a milder growth in the spectral radius
as the order is increased, and the results below confirm this expectation. 

As above, the computational domain $[a,b]$ is divided into non-overlapping uniform DG elements $I^k = [x^k, x^{k+1}], k = 1,\cdots,n$ with element size $H = (b-a)/n$. Each DG element $I^k$ is partitioned into $N$ equidistant subcells with cell size $h = H/N$. Namely, we have $N+1$ degrees of freedom for each DG element. Then the displacement $u(x,t)$ and the velocity $v(x,t)$ in $I^k$ are approximated by the nodal values $U_l^k$ and $V_l^k$ with
\[u(x,t) = \sum_{l=0}^N U_l^k\phi_{p,l}(x),\quad v(x,t) = \sum_{l=0}^N V_l^k\phi_{p,l}(x),\]
respectively. Here, $\phi_{p,l}$ are $p$-th order Galerkin difference basis functions. In this experiment the speed of sound is $c = 1$ and we choose the splitting parameter in the upwind flux to be $\xi = c$. 

Let the vectors $U$ and $V$ contain the nodal values $U_l^k$ and $V_l^k$ 
\[
U = [U_0^1,\cdots,U_N^1,\cdots,U_0^n,\cdots,U_N^n], \quad V = [V_0^1,\cdots,V_N^1,\cdots,V_0^n,\cdots,V_N^n].
\]
Then the semi-discretization (\ref{DG_matrix1})-(\ref{DG_matrix2}) can be written as a system of ordinary differential equations  
\begin{equation*}
\frac{d{Z}}{dt} = \mathcal{L}Z, \quad Z = (U,V)^T.
\end{equation*}
In the experiments in this section, the computational domain is chosen to be $x\in[0,1]$, the number of DG elements is set to be $n = 1$, the degree of approximation space is given by $p = (1,3,5,7,9,11)$, and the number of degrees of freedom in the DG element is given by $N+1 = (31,61,121)$.

\graphicspath{{graphs_spectral_radius/}}
\begin{figure}[h]
	\begin{center}
		\includegraphics[width=0.48\textwidth]{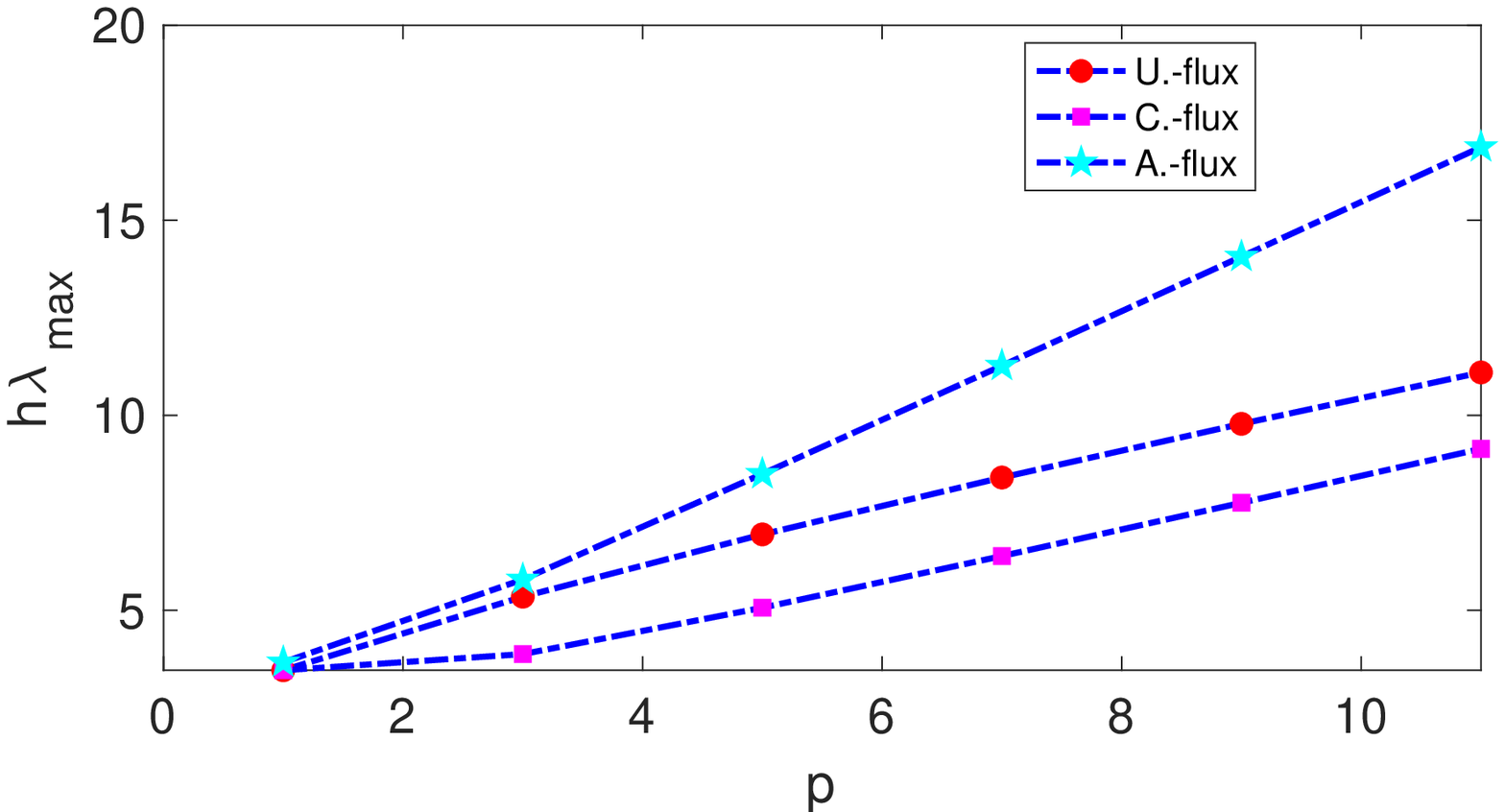}
		\includegraphics[width=0.48\textwidth]{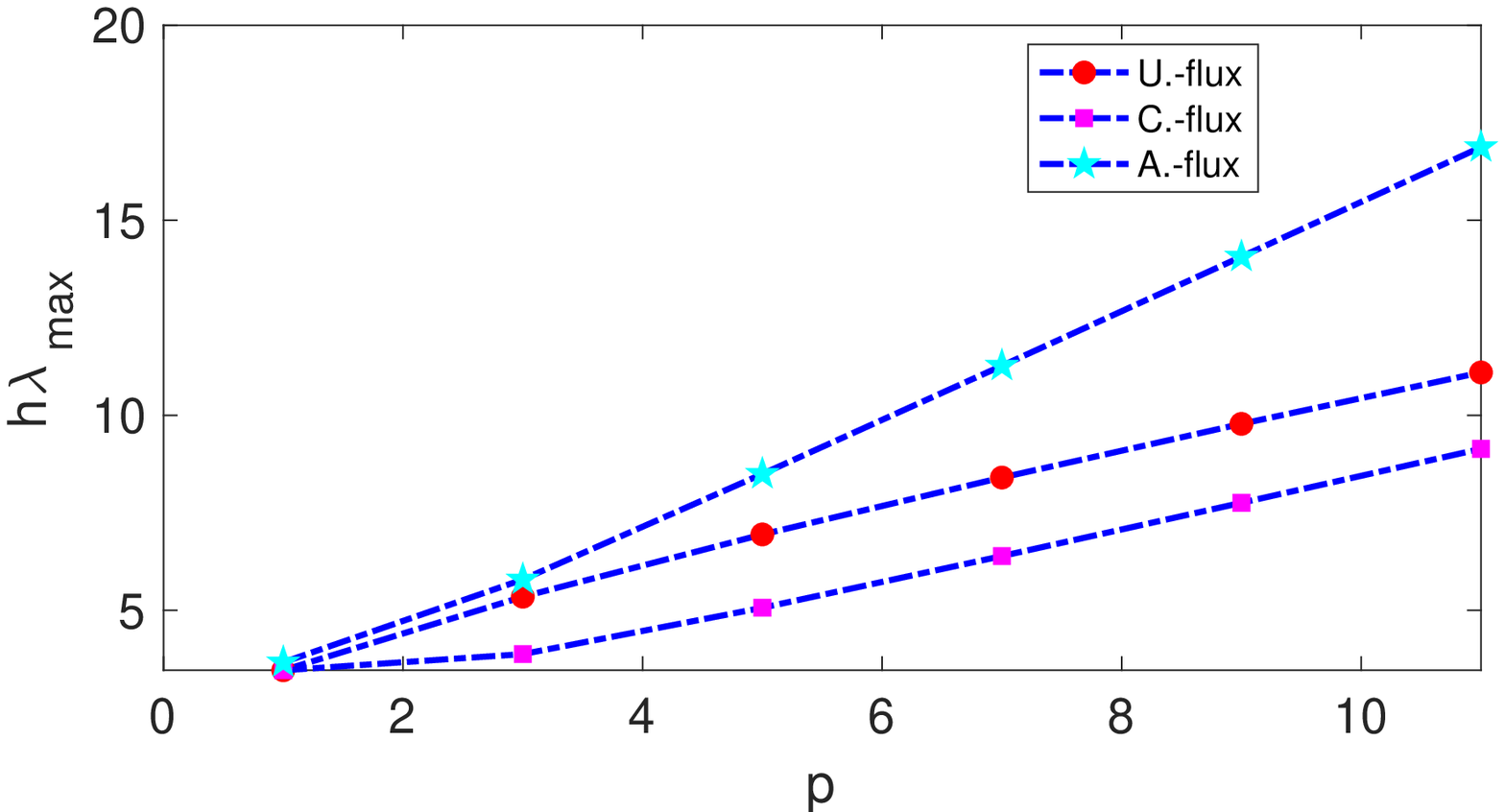}
	\end{center}
	\caption{Spectral radius as a function of the approximation degree $p$ for three different fluxes. The figures are for periodic boundary conditions (left) and Dirichlet and Neumann conditions (right).}\label{spectral_radii_vs_degree_new}
\end{figure}  

Denote the central flux by C.-flux, the alternating flux by A.-flux and the upwind flux by U.-flux. Figure \ref{spectral_radii_vs_degree_new} displays the amplitude of the largest eigenvalue as a function of the degree $p$ for three different values of $N$. The subfigures display results for periodic and non-periodic boundary conditions. We observe that the spectral radius is linearly proportional to the degree $p$ for all three different numerical fluxes. This is in contrast to standard discontinuous finite elements where the spectral radius grows quadratically with $p$. As a consequence the method proposed here can march in time with $p$ times larger time steps. 

\graphicspath{{graphs_convergence_1d/}}

\section{Convergence}\label{convergence}
In this section, we present numerical results to investigate the convergence of our method in both the $L^2$ norm, $\|u-u^h\|_{L^2}$, and the energy norm, $\left(\|\nabla(u - u^h)\|_{L^2}^2 + \|v-v^h\|^2\right)^{1/2}$. We consider both one dimensional and two dimensional problems. As we have discussed in earlier papers \cite{DATH_UP,appelo2020energy,zhang2019energy}, the energy norm is the starting point for theoretically establishing stability and rates of convergence of the energy DG method, and the results there apply directly to the proposed scheme. In all experiments below we use a nodal formulation associated with the basis functions in Section \ref{GD_basis} and march in time by the classic fourth order accurate Runge-Kutta method. We choose the speed of sound to be $c = 1$ and the flux splitting parameter $\xi$ in the upwind scheme to be $c$.

\begin{figure}[htb!]
	\begin{center}
		\includegraphics[width=0.48\textwidth]{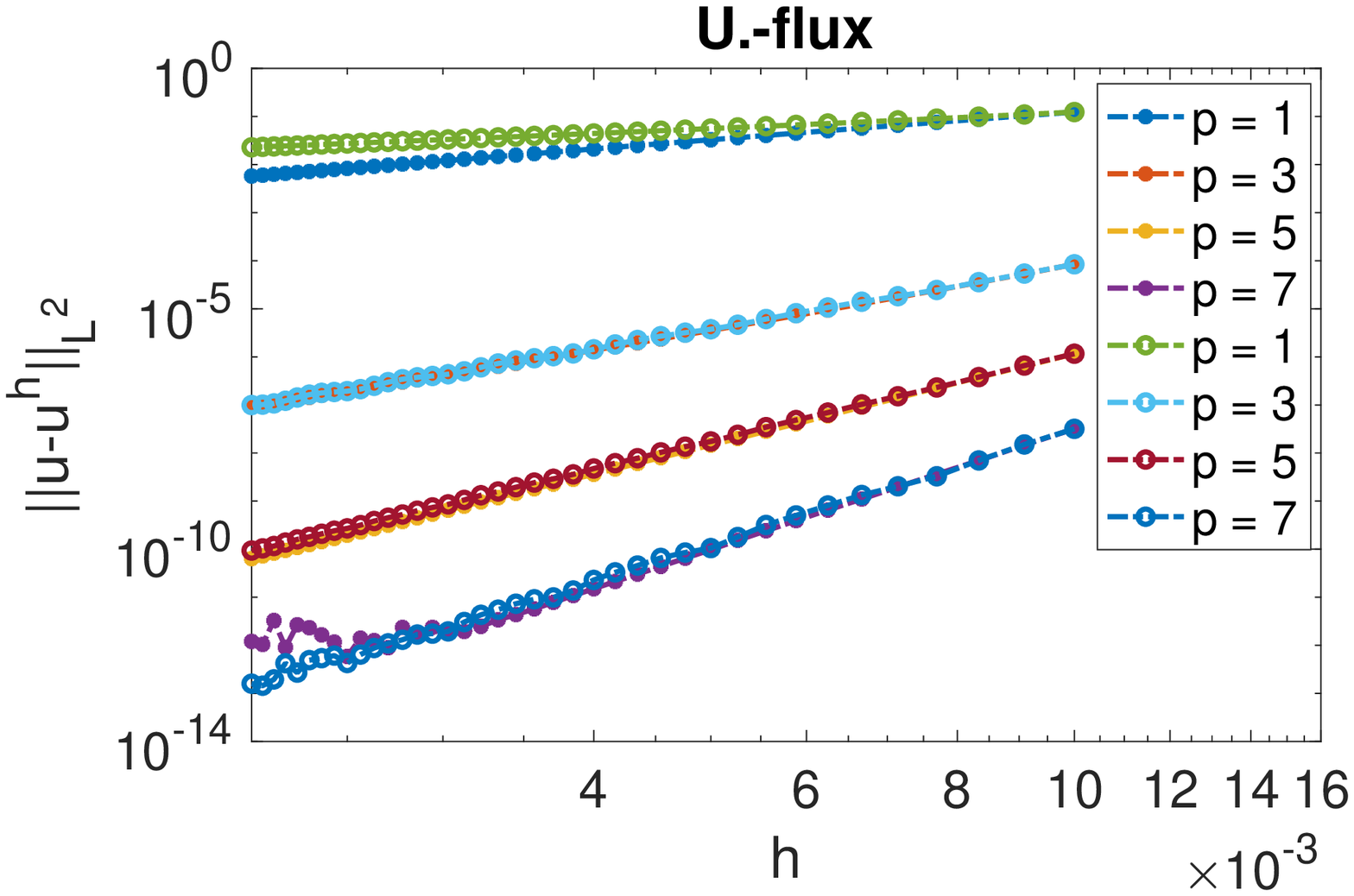}
		\includegraphics[width=0.48\textwidth]{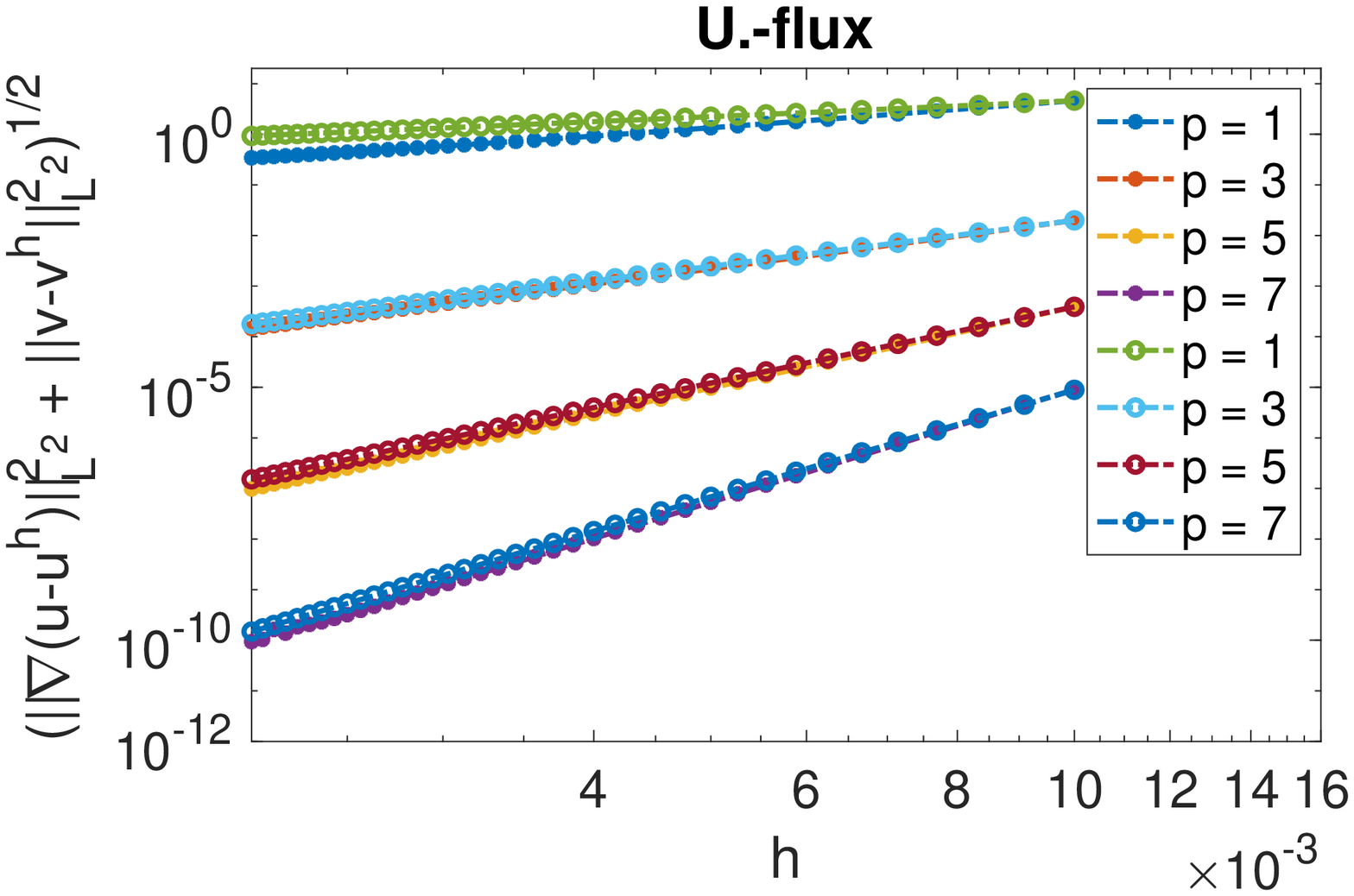}\\
		\includegraphics[width=0.48\textwidth]{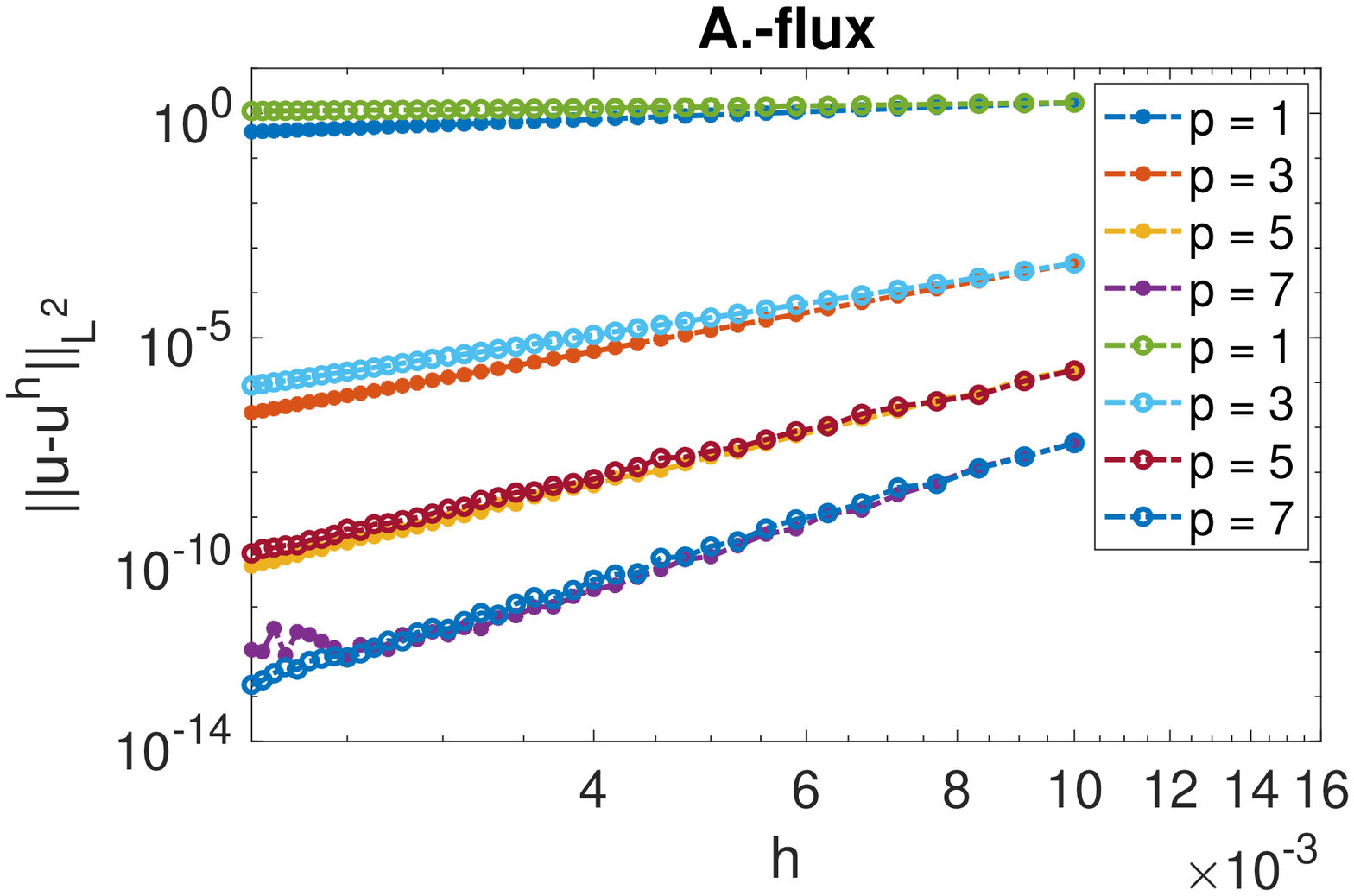}
		\includegraphics[width=0.48\textwidth]{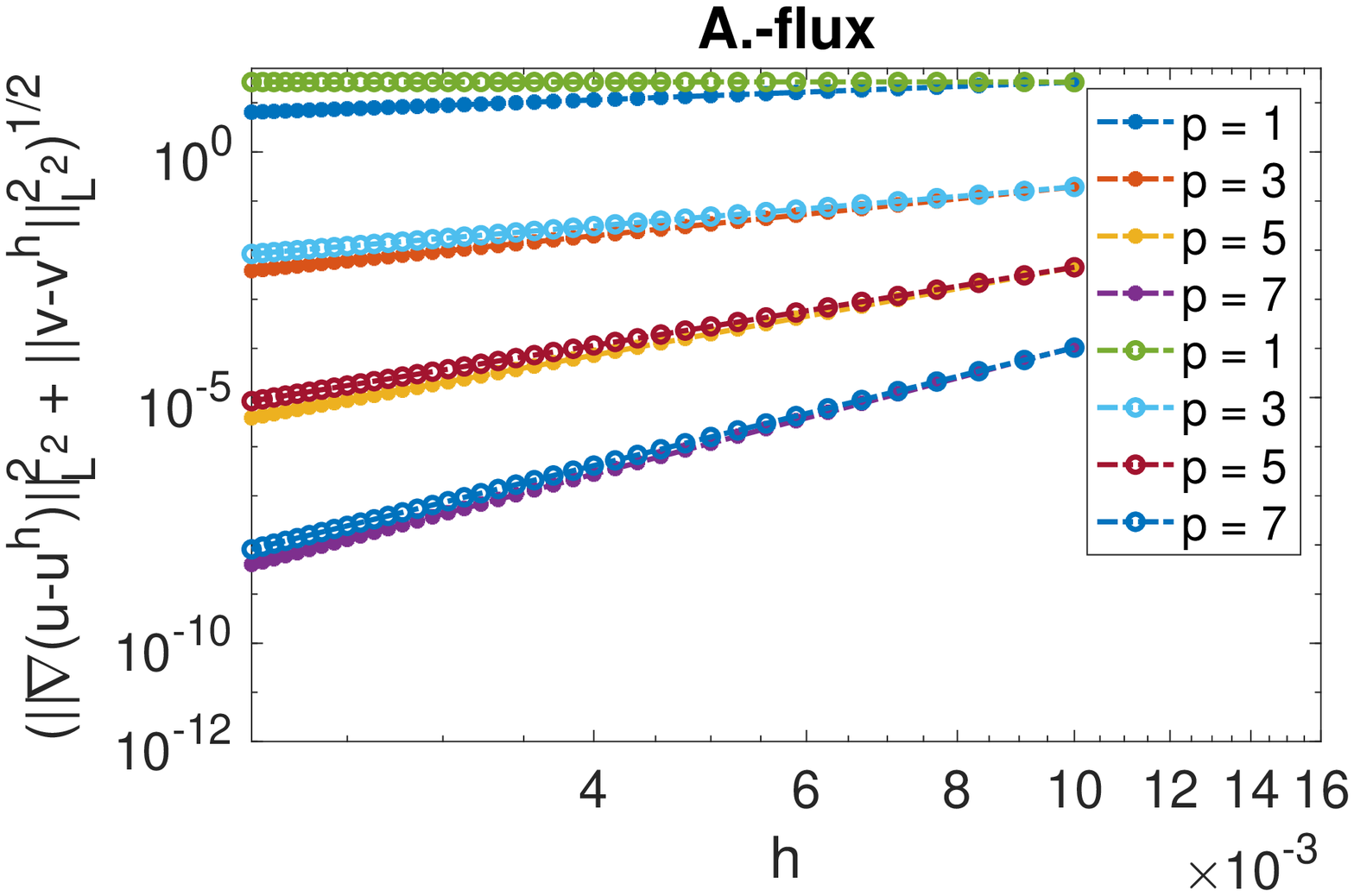}\\
		\includegraphics[width=0.48\textwidth]{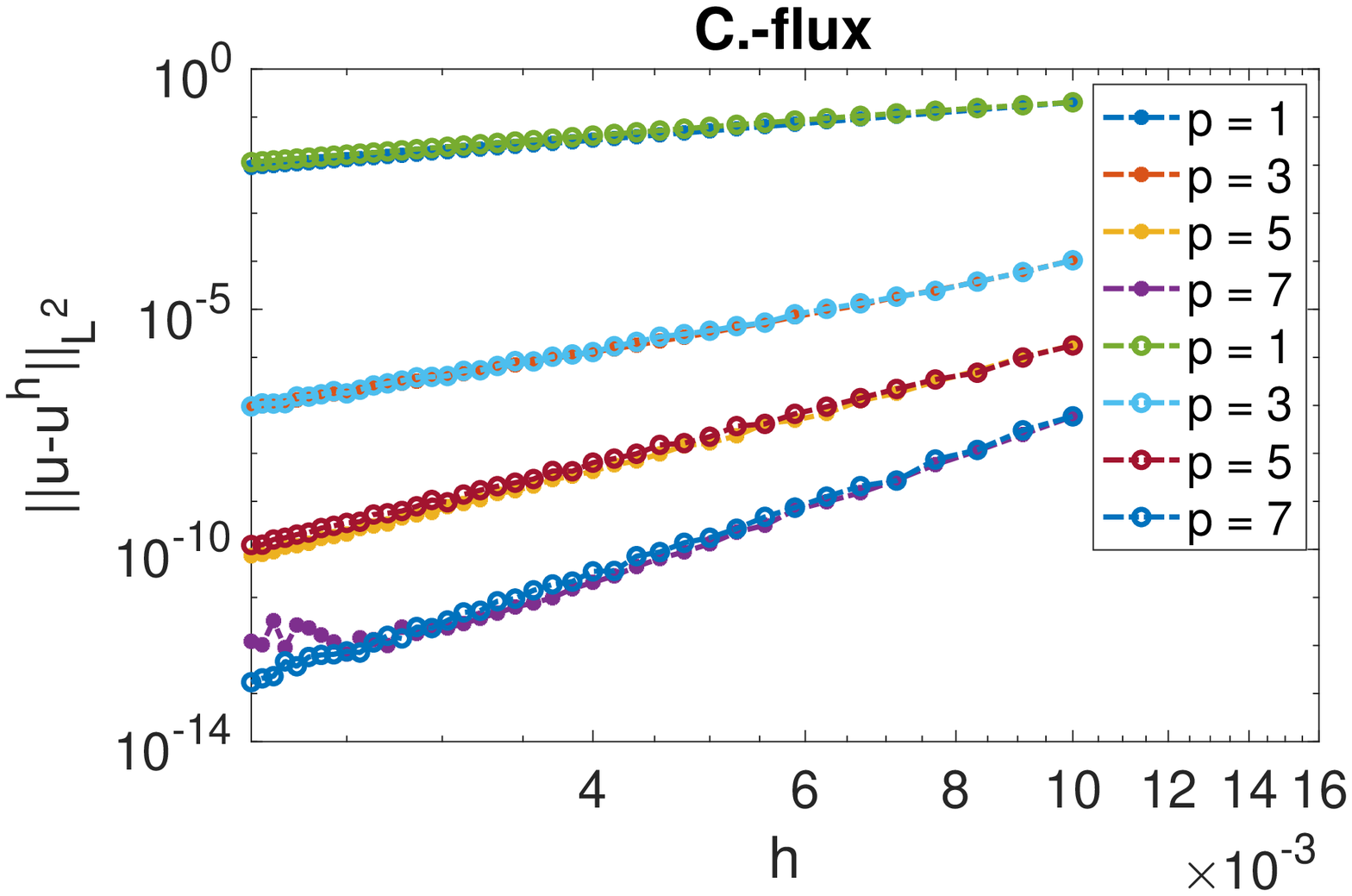}
		\includegraphics[width=0.48\textwidth]{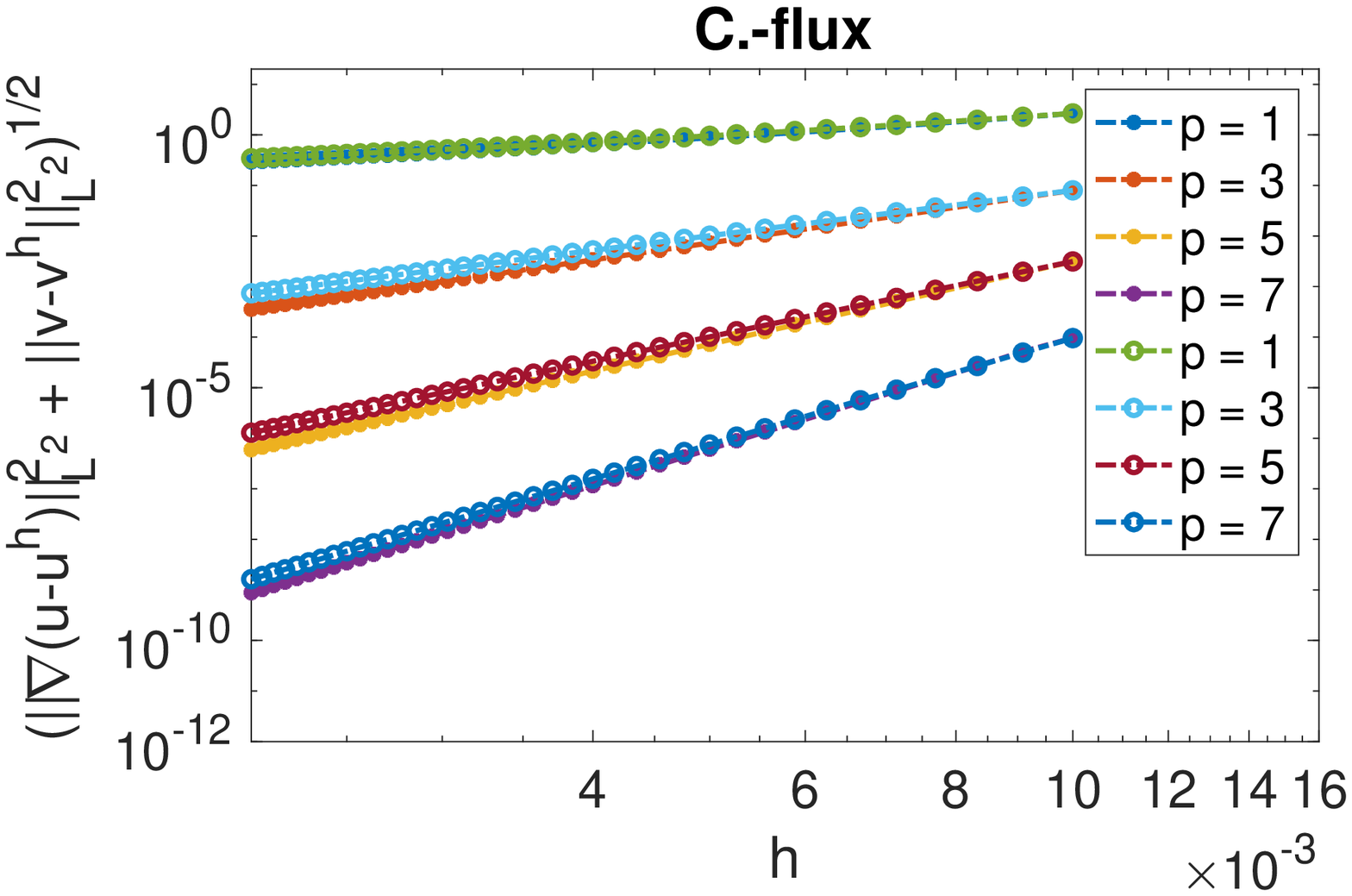}
		\caption{Plots of the $L^2$ norm error $\|u-u^h\|_{L^2}$ and the energy norm error $\left(\|\nabla(u-u^h)\|_{L^2}^2+\|v-v^h\|_{L^2}^2\right)^{1/2}$ as functions of $h$ in one dimension with upwind flux, central flux and alternating flux for periodic boundary conditions. In the legend, $p$ is the degree of the approximation space of $u$ and $v$. The filled circles are the results for the fixed degrees of freedom $N+1 = 11$ and the hollow circles are the results for the fixed number of DG elements $n = 10$.\label{fig:uv_1d_periodic}}
	\end{center}
\end{figure}
\subsection{Periodic Boundary Conditions in One Dimension} \label{numerical_1d_periodic}
To investigate the order of accuracy of our method, we solve
\begin{equation*}
\frac{\partial^2 u}{\partial t^2} = c^2 \frac{\partial^2 u}{\partial x^2}, \quad x\in(0,1),\quad t\geq0,
\end{equation*}
with periodic boundary conditions $u(0,t) = u(1,t)$ and with initial data so that the solution is the exact traveling wave 
\begin{equation}\label{1d_sol_conv}
u(x,t) = \sin(8\pi(x-ct)).
\end{equation}
The discretization is performed on a uniform mesh with DG element vertices $x_i = iH$, $i = 0,\cdots,n$, $H = 1/n$. The vertices of the subcells for the $i$-th DG element are $x_{ij} = x_i + jh$ with $i = 0,\cdots,n-1$, $j = 0,\cdots,N$ and $h = H/N$. The problem is evolved until the final time $T = 1.075$ with the time step size $\Delta t = \mbox{CFL}\times h$, $\mbox{CFL} = 0.075/(2\pi)$ to guarantee the error is dominated by the spatial error. 

We present results for the degree of the approximation space of $u^h, v^h$ being $p = (1,3,5,7)$. The mesh size $h$ is refined by either increasing the number of DG elements $n$ and fixing the degrees of freedom $N+1$ in each DG element or increasing the degrees of freedom $N+1$ in each DG element and fixing the number of DG elements $n$. 

The $L^2$ norm error $\|u-u^h\|_{L^2}$ and the energy norm error $\left(\|\nabla(u-u^h)\|_{L^2}^2+\|v-v^h\|_{L^2}^2\right)^{1/2}$ are plotted against the grid spacing $h$ in Figure \ref{fig:uv_1d_periodic} with the upwind flux, the alternating flux and the central flux, respectively. Linear regression estimates of the rate of convergence can be found in Table \ref{convergence_rate_1d_Nchange} for fixed number of DG elements $n = 10$ and in Table \ref{convergence_rate_1d_Nfixed} for fixed degrees of freedom $N+1 = 11$ in each DG element. Note that we use the same (on element) mesh size $h = 1/n/N$ for these two cases when we generate the results in Figure \ref{fig:uv_1d_periodic}. From Figure \ref{fig:uv_1d_periodic}, we observe that the $L^2$ error $\|u-u^h\|_{L^2}$ becomes oscillatory when it reaches $10^{-13}$ for $p = 7$, thus we only use the first $31$ data points to estimate the convergence rate for this case. For the other cases we use all data to estimate the convergence rate. Generally speaking, the errors are comparable for the two modes of refinement ($N$ fixed and $n$ fixed), with slightly smaller errors if we fix $n$ and refine within each element. From the comparisons of the dispersion errors we do see that refinement within an element is likely to be preferable in most cases. Note that the existing theory for the method in one space dimension only establishes convergence in the energy norm and only proves optimal convergence for the upwind flux. That is, we only have proofs of convergence at the optimal order $p$ in the energy norm for the upwind flux; for the other fluxes the existing theory only guarantees a rate of $p-1$ and in particular it does not guarantee convergence for the conservative fluxes when $p=1$. Specifically we observe the following, which in most cases is better than what we can prove. 

\begin{itemize}
\item[a).] When the number of DG elements $n = 10$ is fixed, from Table \ref{convergence_rate_1d_Nchange}, we observe convergence at rates exceeding $(p+1)$ in the $L^2$ norm for $p = (3,5,7)$ with all three choices for the flux; $p$-th order convergence for the alternating flux and optimal convergence $(p+1)$ for the upwind flux and the central flux when $p = 1$. In the energy norm, we again observe convergence at rates exceeding $p$ for both the upwind flux and the central flux with $p = (1,3,5,7)$, but suboptimal convergence for the alternating flux with $p = (3,5,7)$. We note that we do not expect the asymptotic convergence rates to exceed $p+1$ in $L^2$ and $p$ in the energy norm. However, as we refine within each element, the effect of the $2p$th
order interior formulas is felt, and it is possible to observe convergence at higher rates for some range of
resolutions. This effect is observed for the continuous Galerkin difference methods in \cite{Gdiff}. 
\item[b).] When the degrees of freedom per element $N+1 = 11$ is fixed, from Table \ref{convergence_rate_1d_Nfixed}, we note optimal convergence $(p+1)$ in the $L^2$ norm for $p = (3,5,7)$ with all three fluxes; there is order reduction for $p = 1$. From the energy norm, we have $p$-th order convergence for $p = (1,3,5,7)$ with both the upwind flux and the central flux; for the alternating flux, we observe the suboptimal $(p-1)$ convergence rate for all $p$, and in particular no
convergence when $p=1$. 
\end{itemize}

\begin{table}
	\begin{center}
		\scalebox{0.98}{
	\begin{tabular}{|c|c c c c|c c c c|c c c c|}
		\hline
	    ~ &
		\multicolumn{4}{c|}{U.-flux} &
		\multicolumn{4}{c|}{A.-flux} &
		\multicolumn{4}{c|}{C.-flux} \\
		\hline
		degree $p$ for $u^h$& 1& 3 & 5 & 7 & 1& 3 & 5 & 7& 1& 3 & 5 & 7\\
		\hline
		$L^2$ norm rate &1.97 &4.29 &6.26& 8.23&0.96 &4.88 &6.42& 8.31&1.90 &4.35 &6.43& 8.45\\
		\hline
		energy norm rate &1.65 &3.13 &5.28& 7.33&0.90 &2.50 &4.50& 6.50&1.36 &3.45 &5.47& 7.44\\
		\hline
	\end{tabular}
}
\end{center}
\caption{Linear regression estimates of the convergence rate in the $L^2$ norm $\|u-u^h\|_{L^2}$ and energy norm $\left(\|\nabla(u-u^h)\|_{L^2}^2+\|v-v^h\|_{L^2}^2\right)^{1/2}$ in one dimension with periodic boundary condition for fixed number of DG elements $n = 10$.}\label{convergence_rate_1d_Nchange}
\end{table}

\begin{table}
	\begin{center}
		\scalebox{0.98}{
		\begin{tabular}{|c|c c c c|c c c c|c c c c|}
			\hline
			~ &
			\multicolumn{4}{c|}{U.-flux} &
			\multicolumn{4}{c|}{A.-flux} &
			\multicolumn{4}{c|}{C.-flux} \\
			\hline
			degree $p$ for $u^h$& 1& 3 & 5 & 7 & 1& 3 & 5 & 7& 1& 3 & 5 & 7\\
			\hline
			$L^2$ norm rate &1.06 &4.29 &6.01& 7.95&0.27 &4.00 &5.98& 7.95&1.87 &4.32 &6.09& 8.15\\
			\hline
			energy norm rate &1.02 &3.01 &5.00& 7.01&0.00 &2.00 &4.00& 6.02&1.27 &2.99 &4.98& 7.01\\
			\hline
		\end{tabular}
	}
	\end{center}
	\caption{Linear regression estimates of the convergence rate in the $L^2$ norm $\|u-u^h\|_{L^2}$ and energy norm $\left(\|\nabla(u-u^h)\|_{L^2}^2+\|v-v^h\|_{L^2}^2\right)^{1/2}$ in one dimension with periodic boundary condition for fixed degrees of freedom $N+1 = 11$ in each DG element.}\label{convergence_rate_1d_Nfixed}
\end{table}

\subsection{Dirichlet Boundary Conditions in Two Dimensions}
\graphicspath{{graphs_convergence_2d/}}
\begin{figure}[htb!]
	\begin{center}
		\includegraphics[width=0.48\textwidth]{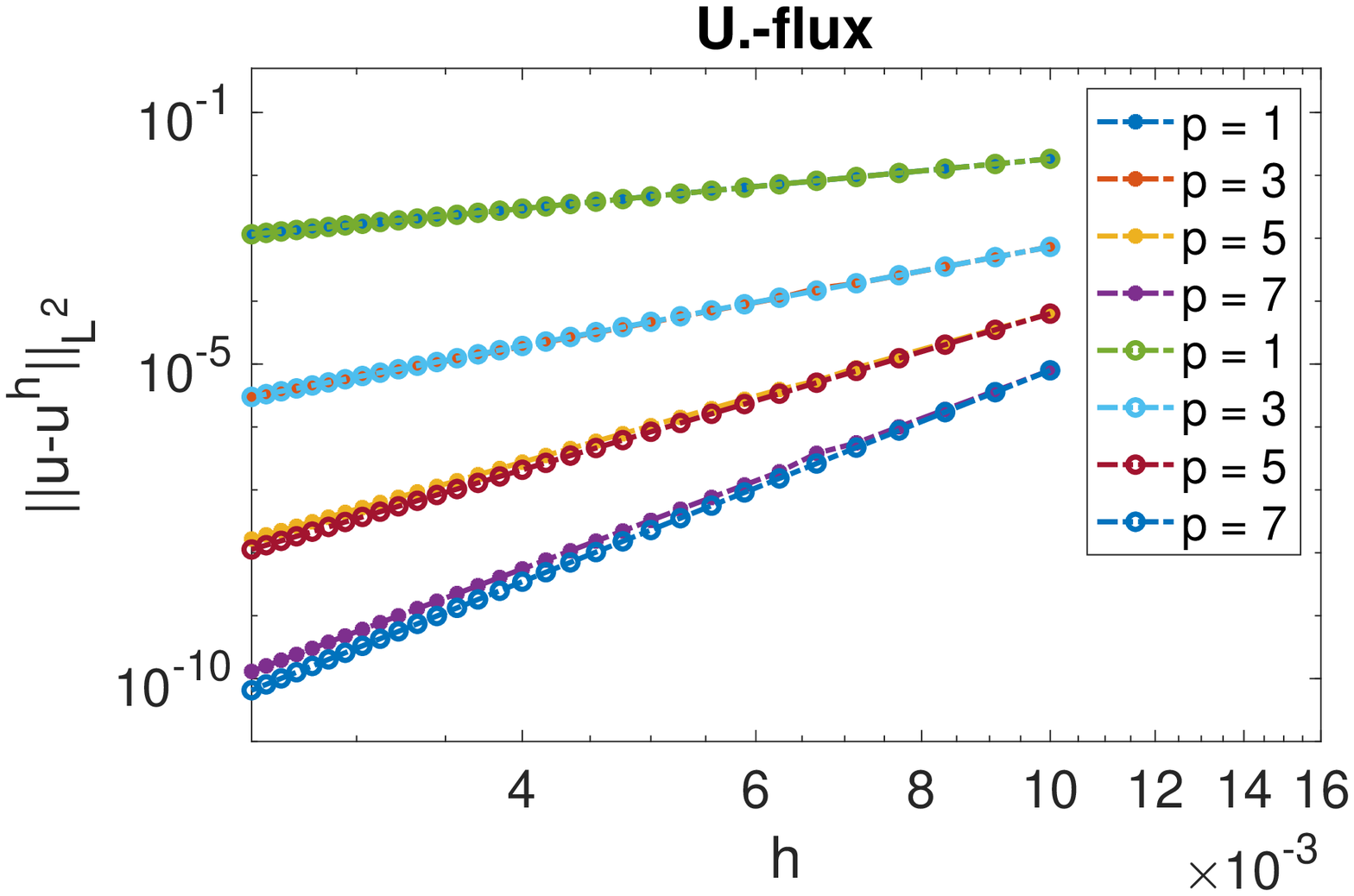}
		\includegraphics[width=0.48\textwidth]{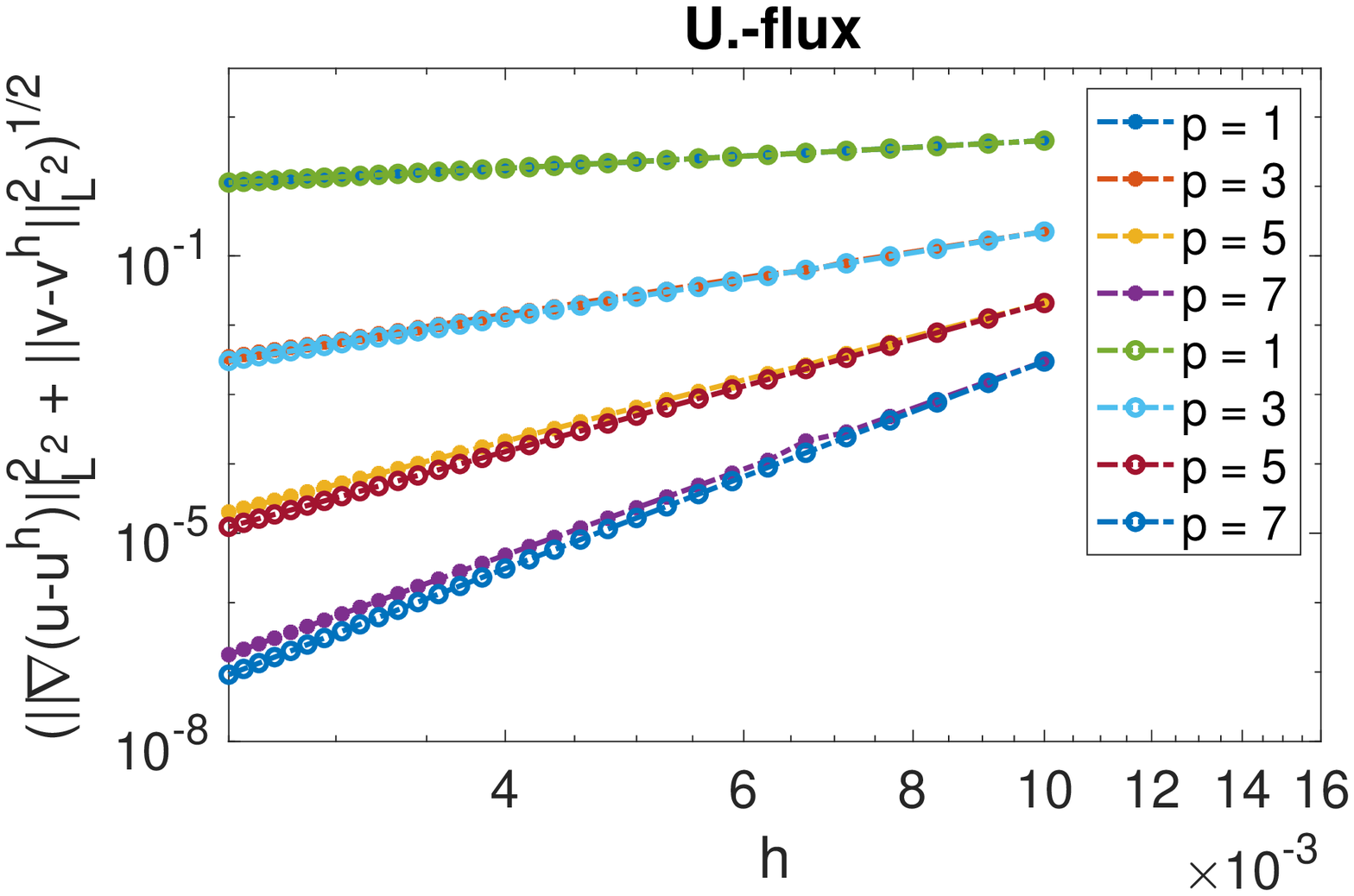}\\
		\includegraphics[width=0.48\textwidth]{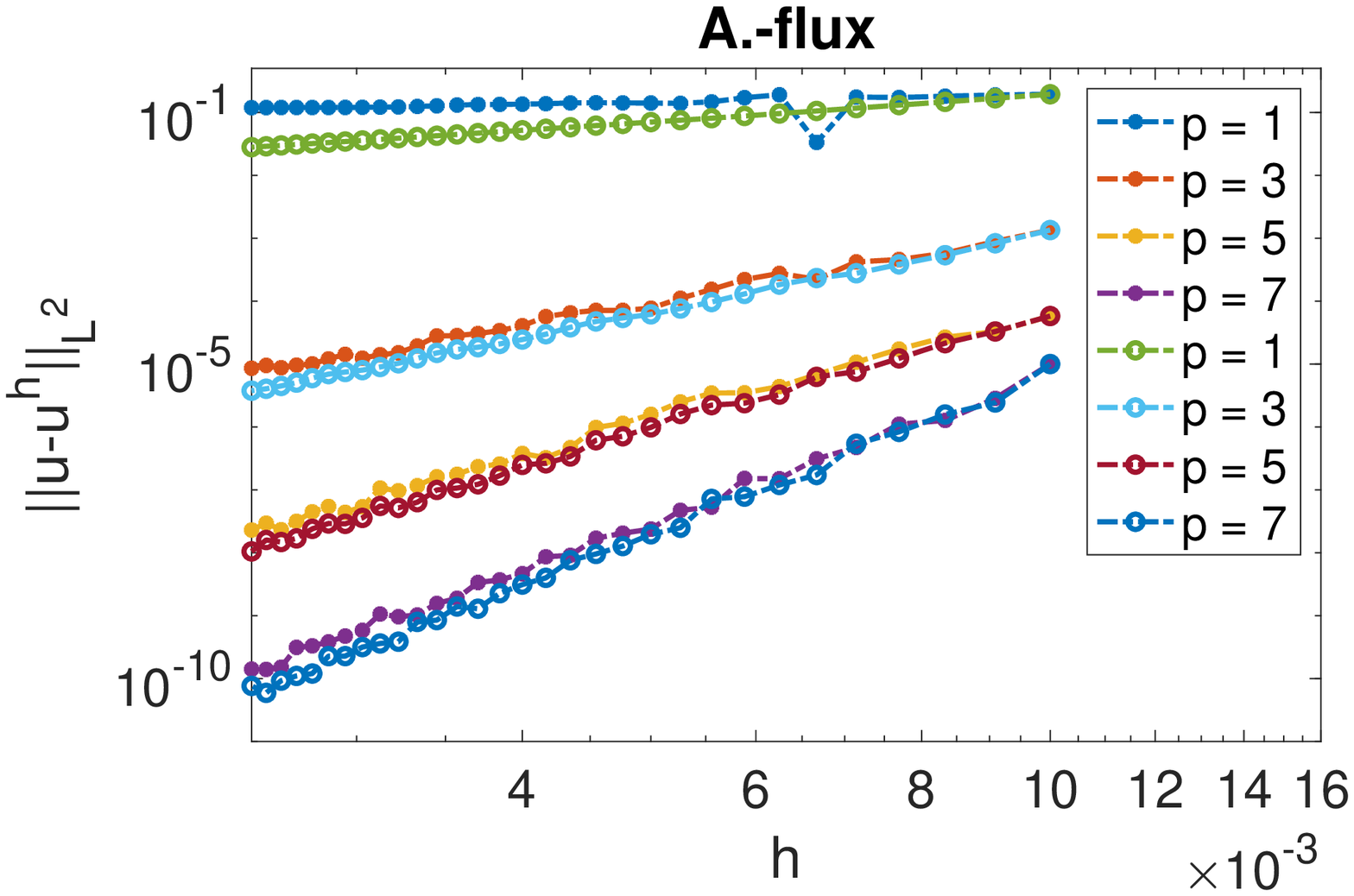}
		\includegraphics[width=0.48\textwidth]{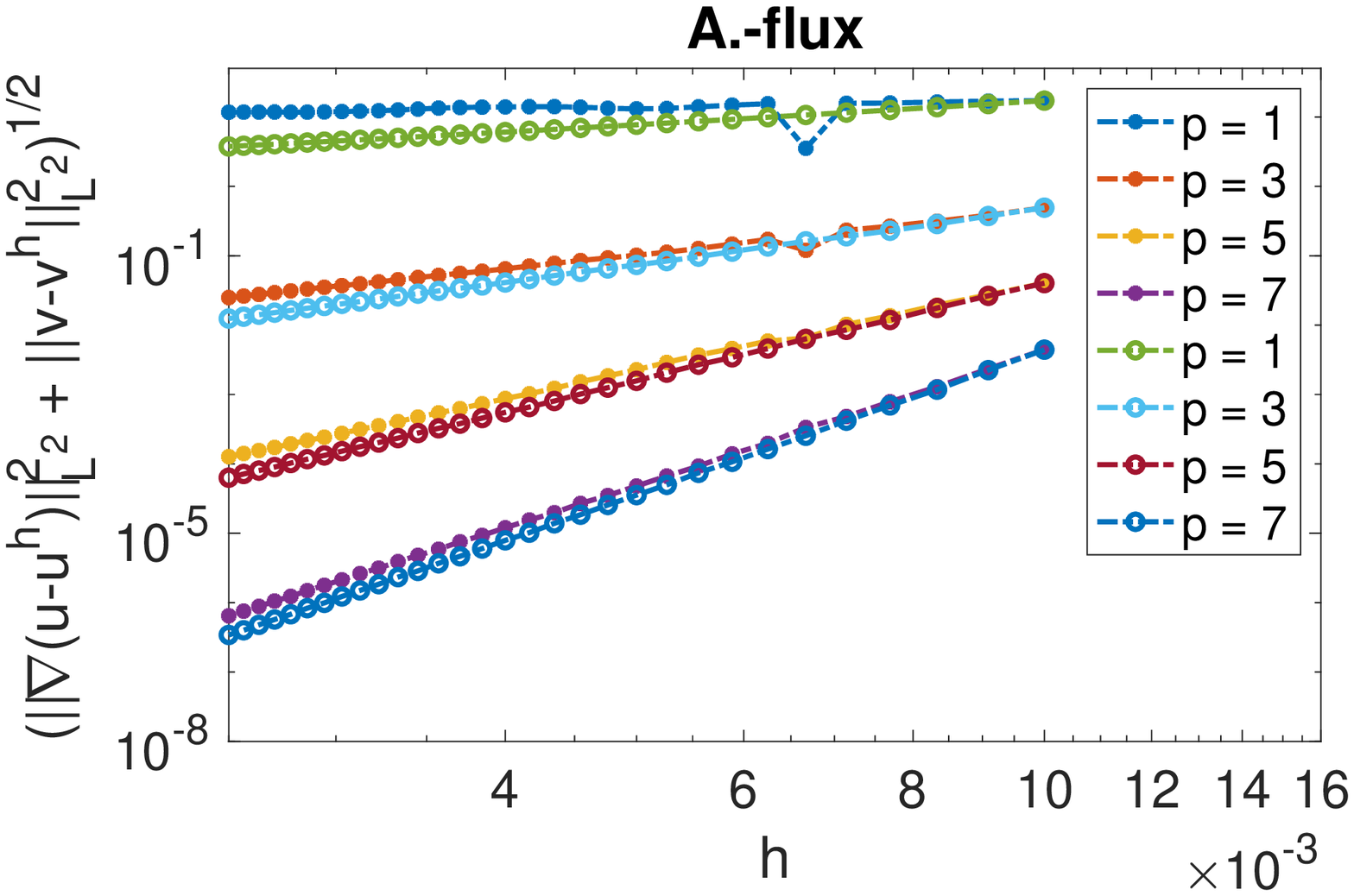}\\
		\includegraphics[width=0.48\textwidth]{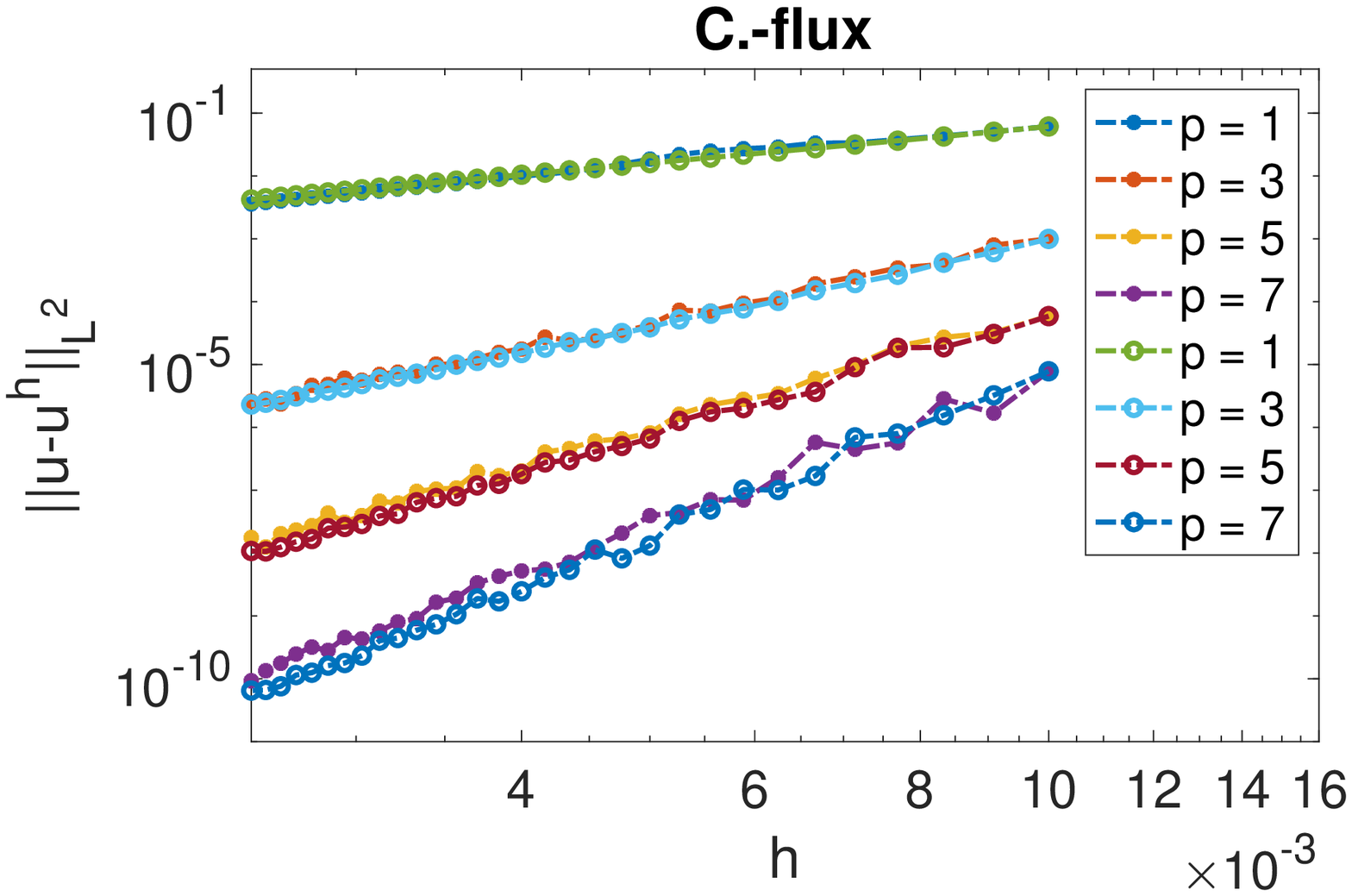}
		\includegraphics[width=0.48\textwidth]{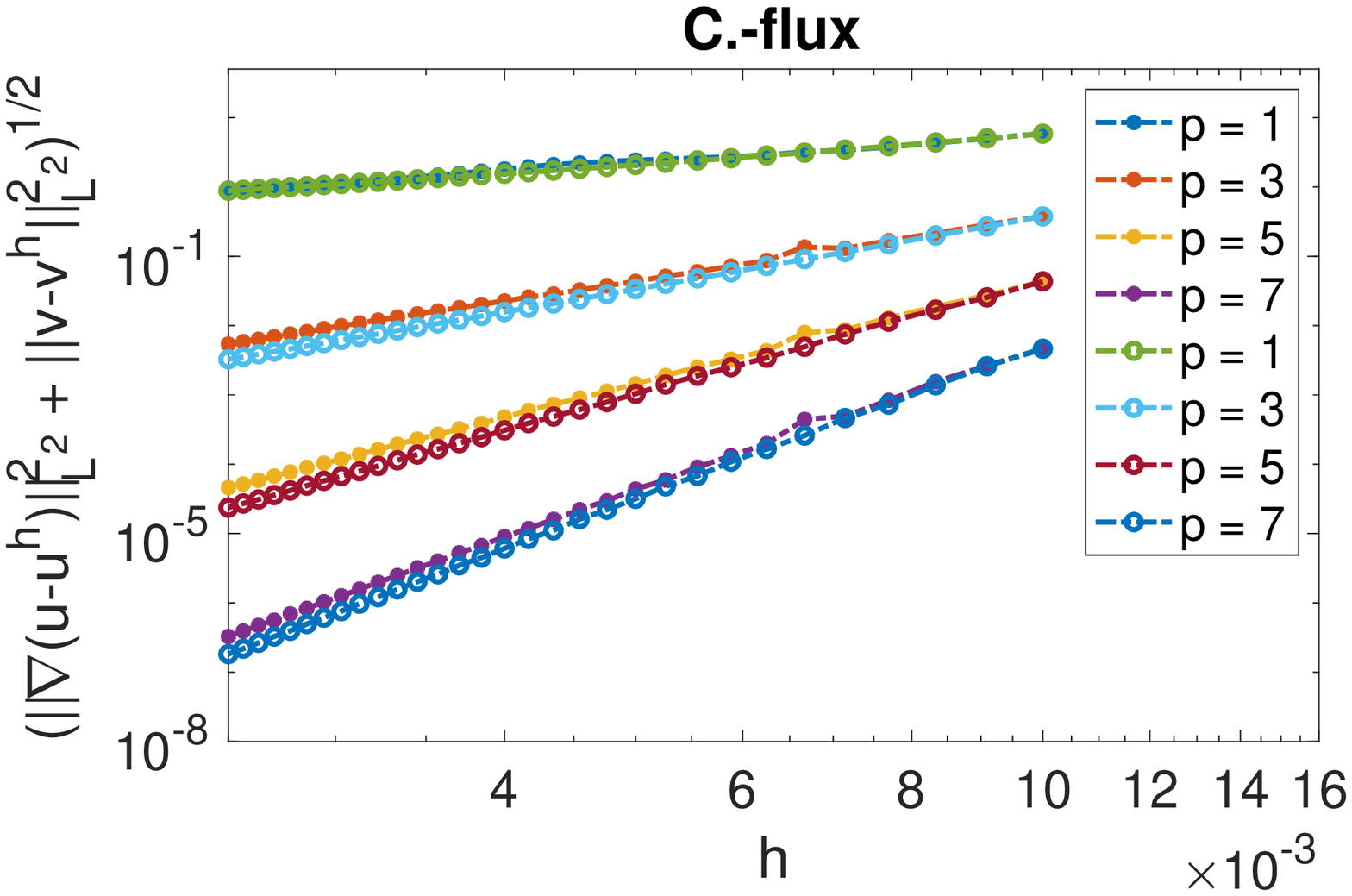}
		\caption{Plots of the $L^2$ norm error $\|u-u^h\|_{L^2}$ and the energy norm error $\left(\|\nabla(u-u^h)\|_{L^2}^2+\|v-v^h\|_{L^2}^2\right)^{1/2}$ as functions of $h$ in two dimensions with upwind flux, central flux and alternating flux for the Dirichlet boundary condition. In the legend, $p$ is the degree of the approximation space of $u$ and $v$. The filled circles are the results for the fixed degrees of freedom $(N+1)^2 = 11^2$ and the hollow circles are the results for the fixed number of DG elements $n^2 = 10^2$.}\label{fig:uv_2d_dirichlet}
	\end{center}
\end{figure}
In this section, we examine the rate of convergence for our scheme applied to the following two dimensional problem  
\begin{equation*}
\frac{\partial^2 u}{\partial t^2} = \left(\frac{\partial^2 u}{\partial x^2} + \frac{\partial^2 u}{\partial y^2}\right), \ \ \ (x,y)\in(0,1)\times(0,1),\ \ t\geq 0,
\end{equation*}
with initial data and boundary conditions chosen such that the exact solution is
\begin{equation*}
u(x,y,t) = \cos(15\sqrt{2}\pi ct)\sin(15\pi x)\sin(15\pi y).
\end{equation*}

The computational domain is discretized into Cartesian DG elements whose vertices are $x_i = iH$, $y_j = jH$, $i,j = 0,1,\cdots n$, $H = 1/n$. The vertices of subcells for the $ij$-th DG element are $x_{ik} = x_i + kh$ and $y_{jl} = y_j + lh$ with $i,j = 0,\cdots,n-1$ and $k,l = 0,\cdots,N$,  $h = H/N$. The problem is evolved until final time $T = 0.2$ and the time step size is $\Delta t = \mbox{CFL} h$, with $\mbox{CFL} = 0.075/(2\pi)$ to guarantee the temporal error is dominated by the spatial error. 

As in the one dimensional test in Section \ref{numerical_1d_periodic}, we consider three different numerical fluxes: the upwind flux, the alternating flux, and the central flux. We again use two different ways to refine the spatial mesh size: one is to fix the degrees of freedom in each element, $(N+1)^2 = 11^2$ and change the number of DG elements $n^2$. The other is to fix the number of DG elements $n^2 = 10^2$ and change the degrees of freedom in each DG element. Again, we have the same (on element) mesh size $h = 1/n/N$ for these two cases.  

\begin{table}
	\begin{center}
		\scalebox{0.98}{
			\begin{tabular}{|c|c c c c|c c c c|c c c c|}
				\hline
				~ &
				\multicolumn{4}{c|}{U.-flux} &
				\multicolumn{4}{c|}{A.-flux} &
				\multicolumn{4}{c|}{C.-flux} \\
				\hline
				degree $p$ for $u^h$& 1& 3 & 5 & 7 & 1& 3 & 5 & 7& 1& 3 & 5 & 7\\
				\hline
				$L^2$ norm rate &1.99 &3.96 &6.23& 8.46&1.34 &4.18 &6.24& 8.47&1.92 &4.33 &6.37& 8.58\\
				\hline
				energy norm rate &1.00 &3.09 &5.35& 7.51&1.06 &2.63 &4.68& 6.80&1.34 &3.42 &5.45& 7.39\\
				\hline
			\end{tabular}
		}
	\end{center}
	\caption{Linear regression estimates of the convergence rate in the $L^2$ norm $\|u-u^h\|_{L^2}$ and energy norm $\left(\|\nabla(u-u^h)\|_{L^2}^2+\|v-v^h\|_{L^2}^2\right)^{1/2}$ in two dimensions with Dirichlet boundary condition for fixed degrees of freedom per element $N+1 = 11$.}\label{convergence_rate_2d_Nchange}
\end{table}

\begin{table}
	\begin{center}
		\scalebox{0.98}{
			\begin{tabular}{|c|c c c c|c c c c|c c c c|}
				\hline
				~ &
				\multicolumn{4}{c|}{U.-flux} &
				\multicolumn{4}{c|}{A.-flux} &
				\multicolumn{4}{c|}{C.-flux} \\
				\hline
				degree $p$ for $u^h$& 1& 3 & 5 & 7 & 1& 3 & 5 & 7& 1& 3 & 5 & 7\\
				\hline
				$L^2$ norm rate &1.99 &4.00 &5.95& 7.97&0.39 &3.74 &5.85& 7.86&2.14 &4.33 &6.14& 7.96\\
				\hline
				energy norm rate &1.00 &2.99 &5.00& 7.04&0.27 &2.07 &4.18& 6.34&1.28 &3.06 &5.00& 6.99\\
				\hline
			\end{tabular}
		}
	\end{center}
	\caption{Linear regression estimates of the convergence rate in the $L^2$ norm $\|u-u^h\|_{L^2}$ and energy norm $\left(\|\nabla(u-u^h)\|_{L^2}^2+\|v-v^h\|_{L^2}^2\right)^{1/2}$ in two dimensions with Dirichlet boundary condition for fixed number of DG elements $n = 10$.}\label{convergence_rate_2d_Nfixed}
\end{table}

The $L^2$ norm error $\|u-u^h\|_{L^2}$ and the energy norm error $\left(\|\nabla(u-u^h)\|_{L^2}^2+\|v-v^h\|_{L^2}^2\right)^{1/2}$ are presented in Figure \ref{fig:uv_2d_dirichlet} for the upwind flux, the central flux and the alternating flux, respectively. The corresponding convergence rates from linear regression are shown in Table \ref{convergence_rate_2d_Nchange} for a fixed number of DG elements $n^2 = 10^2$ and in Table \ref{convergence_rate_2d_Nfixed} for fixed degrees of freedom $(N+1)^2 = 11^2$ on each DG element. We use all data to generate the convergence rate here, but for the alternating flux with fixed degrees of freedom $(N+1)^2 = 11^2$, we use the data from the $25$ coarsest grids which excludes the outliers where the error is very small. Generally speaking, the results are similar to the one dimensional results in Section \ref{numerical_1d_periodic}. We observe optimal convergence when $p = (3,5,7)$ for all cases. However, we observe a rate of convergence $p+1$ in the $L^2$ error norm for the upwind flux when $p = 1$ with fixed degrees of freedom $(N+1)^2 = 11^2$ on each DG element. This is slightly better than the corresponding one dimensional result.

\section{Problems with Variable Coefficients}\label{variable}
As was demonstrated in Section \ref{GD_cost} when applied to constant coefficient problems on Cartesian meshes the proposed method has the same complexity as a traditional finite difference method. Unfortunately the simultaneous diagonalization cannot be expected to work when the speed of sound varies in space or for a constant coefficient problem on non-Cartesian meshes (in the latter case the transformation from a physical element to the reference element will result in a variable coefficient problem). For such variable coefficient problems we will stay with the standard Galerkin difference basis and invert the matrices on the left of the element-wise equations (\ref{DG_matrix1}) and (\ref{DG_matrix2}) using the preconditioned conjugate gradient (pcg) method. In a first experiment we demonstrate that the number of pcg iterations needed in each time step is small. In a second experiment we demonstrate the ability of the method to compute the solution to a more complex application-type problem.

\subsection{Efficiency of PCG for Inverting Mass Matrices}\label{iterative_variable}
Here, we consider the second order wave equations with a variable coefficient $c^2(x,y)$ in two dimensions as follows
\begin{equation}\label{variable_problem}
\frac{\partial^2 u}{\partial t^2} = \nabla\cdot(c^2(x,y)\nabla u) + f(x,y), \quad (x,y)\in[0,1]\times[0,1],\quad t\geq0,
\end{equation}
where $c^2(x,y) = 1+x^2+y^2$. The initial conditions and the external forcing function are determined by the manufactured solution
\[
u(x,y,t) = \sin(8\sqrt{2}\pi t)\sin(8\pi x)\sin(8\pi y).
\]
We impose periodic boundary conditions. 

The key point for the success of the new basis functions proposed in Section \ref{GD_new_basis} is the mass matrix and stiffness matrix in  d-dimensions are constructed by the tensor products of the corresponding mass and stiffness matrices in one dimension. From the scheme (\ref{DG1})-(\ref{DG2}), though the mass matrix for the problem (\ref{variable_problem}) still keeps the tensor product form, the elements of the stiffness matrix $S_{c^2(x,y)}$ are derived from
\[
\int_{\Omega_k} c^2(x,y)\nabla\varphi_u\cdot\frac{\partial \nabla u}{\partial t} d\Omega^k,
\]
which does not have a tensor product form for general $c^2(x,y)$. In \cite{chan2017weight,chan2017weight2}, the authors proposed weight-adjusted discontinuous Galerkin (WADG) method to handle the variable coefficient in the mass matrix. The idea there is to replace the weighted $L^2$ inner product with a weight-adjusted inner product. Unfortunately this approach is not applicable here.

As an alternative for the variable coefficient $c^2(x,y)$ problem, we compute the time derivatives in the scheme (\ref{DG_matrix1})--(\ref{DG_matrix2}) iteratively by the preconditioned conjugate gradient. As a preconditioner we use the zero fill-in incomplete Cholesky factorization of $S_{c^2(x,y)}$ and $M$ as the preconditioning matrix of the system (\ref{DG_matrix1}) and (\ref{DG_matrix2}), respectively. In particular we use the Fortran subroutines of Jones and Plassmann \ci{JonesPlassmannIC}. The mesh, time stepping and other parameters are the same as in the two dimensional example above. The degree of the approximation space of $u$ and $v$ is chosen to be $p = 3$.

\begin{table}
	\begin{center}
		\scalebox{0.98}{
	\begin{tabular}{|c|c| c| c| c| c|}
		\hline
		\multicolumn{2}{|c|}{Degrees of freedom $(N+1)^2$ in one DG element} &121 & 441 & 1681 & 6561\\
	    \hline
		\multirow{2}{*}{---}&$L^2$ error in $u$  & 3.14e-2& 1.58e-3 & 1.07e-4 & 6.89e-6\\
		\cline{2-6}
		&convergence rate  & --& 4.31& 3.88 & 3.96\\
		\hline
		\multicolumn{2}{|c|}{relative tolerance in PCG iterative method}&$10^{-3}$ & $10^{-3}/16$ & $10^{-3}/16^2$ & $10^{-3}/16^3$\\
		\hline
		PCG for ${dU}/{dt}$ 
	    &average number of iterations &1.58 &2.13 &3.18& 5.03\\
		\hline
		PCG for ${dV}/{dt}$ 
		&average number of iterations &1 &1 &1 & 1\\
		\hline
	\end{tabular}
}
\end{center}
\caption{The average number of iterations for preconditioned conjugate gradient methods in solving $\frac{dU}{dt}$ and $\frac{dV}{dt}$ for variable $c^2(x,y)$ in two dimensions with fixed number of DG elements $n^2 = 4$.}\label{iteration_number}
\end{table}

\graphicspath{{graphs_sofar/}}
\begin{figure}[htb!]
        \begin{center}
                \includegraphics[width=0.9\textwidth]{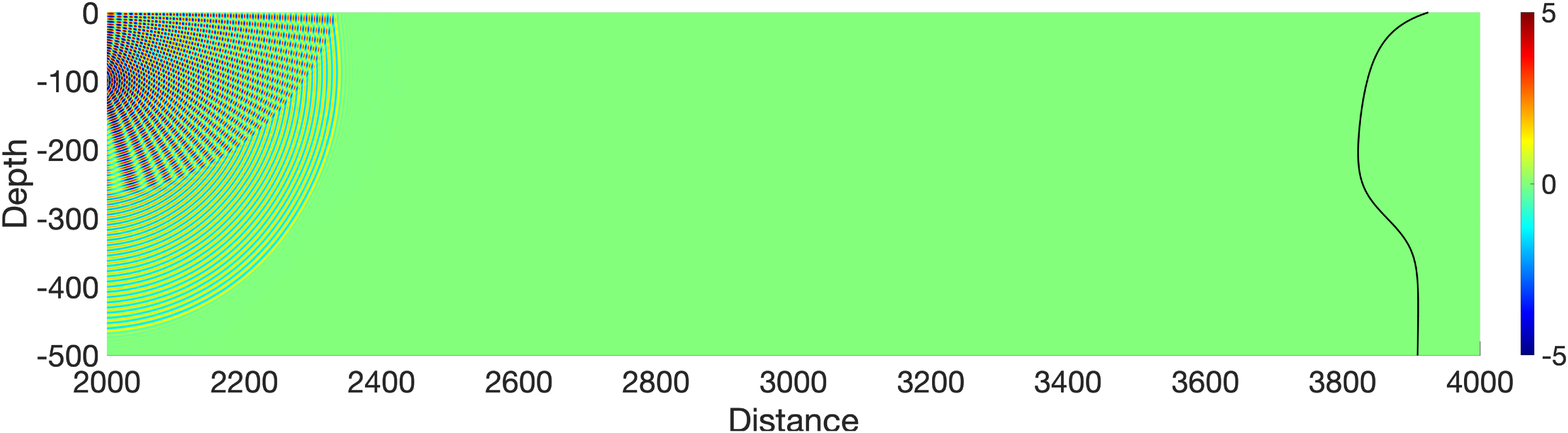}
                \includegraphics[width=0.90\textwidth]{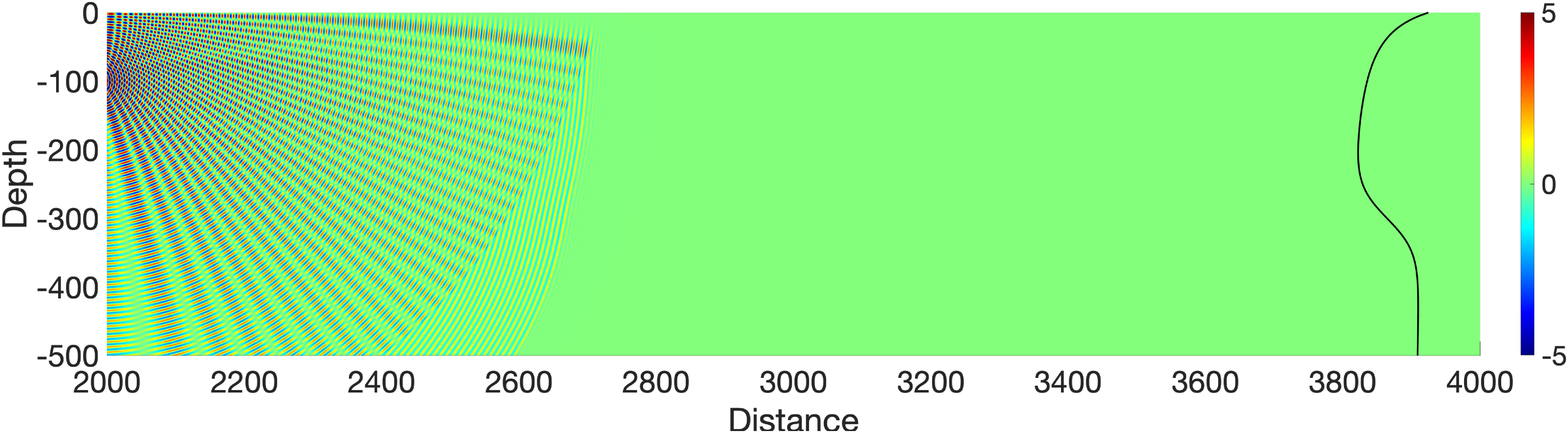}
                \includegraphics[width=0.90\textwidth]{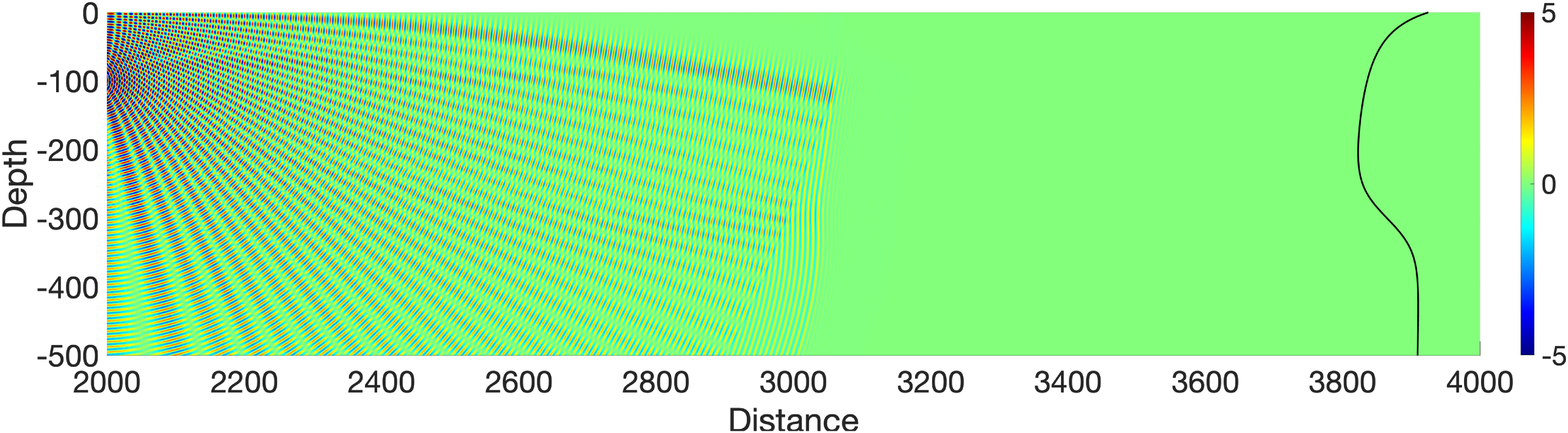}
                \includegraphics[width=0.90\textwidth]{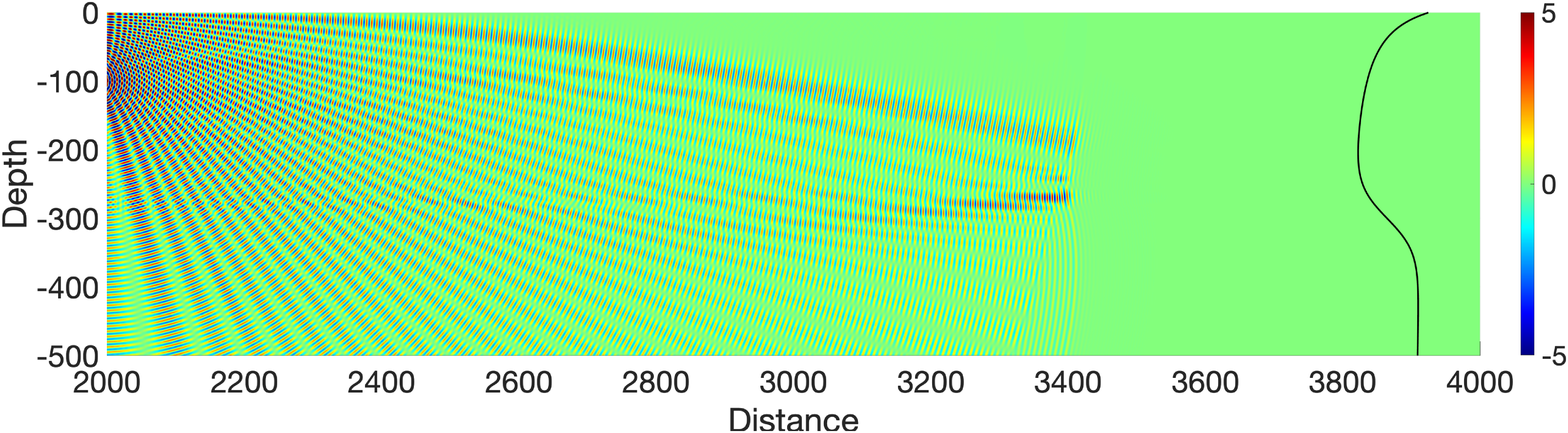}
                \includegraphics[width=0.90\textwidth]{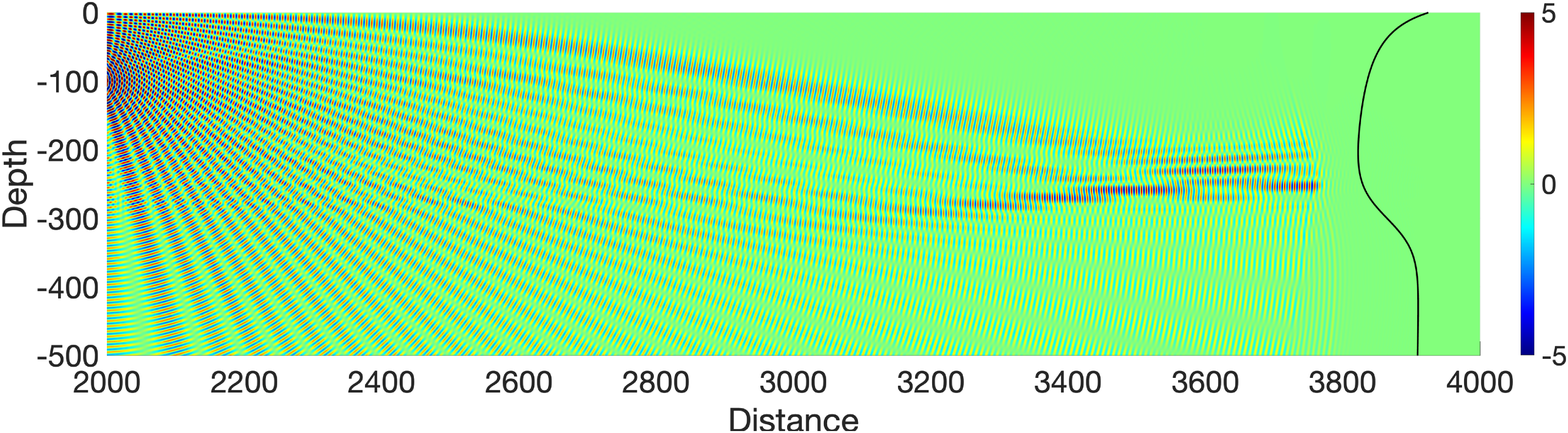}
        \end{center}
        \caption{Time snapshots at times 0.25s, 0.5s, 0.75s, 1.0s  and 1.25s. As the waves progress from left to right the guiding due to the variable speed of sound becomes increasingly pronounced.  The black inlay is a sketch of the sound speed profile. \label{fig:sofar}}
\end{figure}

In Table \ref{iteration_number}, we show the $L^2$ errors in $u$, the corresponding convergence rate and the average number of the iterations for solving systems of $\frac{dU}{dt}$ and $\frac{dV}{dt}$. We observe the $4$-th order convergence in the $L^2$ norm for $u$ which is the same with results for the constant variable $c^2$. In addition, the average number of iterations is comparably small for all different degrees of freedom $(N+1)^2$ on each DG element. The average number of the iterations for $\frac{dU}{dt}$ is less than $6$ and the average number of the iterations for $\frac{dV}{dt}$ is always $1$. From these results (which are largely representative for many experiments that we have conducted), we see that the iterative methods for solving the scheme (\ref{DG_matrix1})--(\ref{DG_matrix2}) with variable stiffness matrix or variable mass matrix work well and are not cost-prohibitive. However the observed increase in the number of iterations with the number of degrees-of-freedom per element, which is
expected for the sort of preconditioner we are using, suggests that this number should be held fixed and the number of DG elements increased as we
refine the mesh.

\subsection{Propagation of waves in an ocean channel}\label{nice_phenomenon}
As a final and more applied example we consider the propagation of sound from a point source 100 meters below the surface in the ocean. The point source has a $\sin (400\pi t)$ sinusodial time dependence.  The speed of sound in water is taken to be depth dependent with the formula for the speed of sound in water being 
$$
c = 1450+50 \left( \frac{100}{D+40} + \tanh \left(\frac{D-300}{50}\right) \right).
$$
Here $D$ is the depth (with a positive value) below the surface. As the sound speed profile has an inflection point the waves will be guided as they propagate in the direction parallel to the surface. Depending on the ``strength'' of the guiding the effect will become visible at a few or many wavelengths. The profile we use here is relatively weak and thus it is important to have a numerical method that is able to accurately propagate waves with minimal error over long distances.

Here the computational domain is taken to be $(x,y) \in [0m, 4000m] \times [-2000m,0m]$ and is discretized by $100 \times 50$ elements each with $41 \times 41$ grid points. We take $p=3$ and use the same Runge-Kutta method as above. In Figure \ref{fig:sofar} we display the solution $u$ at different snapshots in time. As the wave fronts evolve to the right in the domain they concentrate in side the minima of the sound speed profile illustrating that the method is able to capture these phenomena well.     

We note that for the purely depth-dependent sound speed profile considered here the diagonalization method used for constant coefficient problems
can be used. However, to illustrate its utility for a more complex problem we used the iterative method in this case also.

\section{Summary}\label{conclusion}
In conclusion, we have demonstrated the energy-based DG method with Galerkin difference basis functions for second-order wave equations. In particular:

\begin{itemize}
\item[a)] We derived a new basis by simultaneous diagonalization of the mass and the stiffness matrices from the Galerkin difference basis functions. The new basis reduces the computational cost of evolving the solution from superlinear complexity with respect to degrees of freedom to optimal linear complexity.
\item[b)] Using Bloch wave analysis we computed the dispersive and dissipative properties of the method.
\item[c)] By numerical experiments we showed that the spectral radius of the semi-discretization of our scheme is linearly proportional to the degree of the approximation space $p$. This translates to the ability to march in time using $p$ times  larger time steps compared with traditional element based methods such as spectral-, continuous- and discontinuous finite elements.
\item[d)] Optimal convergence was observed for problems in both one and two dimensions for all numerical fluxes when the degree $p\geq3$. The results apply both when
the number of points-per-element is refined with the number of elements fixed and when the number of points-per-element is fixed and the number of elements is
increased.
\item[e)] We illustrated that the method is not dramatically slower for variable coefficient problems if the mass matrices are inverted using the preconditioned conjugated gradient method. In this case fixing the number of points-per-element while increasig the number of elements may be the most efficient refinement
strategy. 
\end{itemize}

\section*{Conflict of interest}

The authors declare that they have no conflict of interest.

\bibliography{GD_newbasis}

\begin{thebibliography}{10}
\providecommand{\url}[1]{{#1}}
\providecommand{\urlprefix}{URL }
\expandafter\ifx\csname urlstyle\endcsname\relax
  \providecommand{\doi}[1]{DOI~\discretionary{}{}{}#1}\else
  \providecommand{\doi}{DOI~\discretionary{}{}{}\begingroup
  \urlstyle{rm}\Url}\fi

\bibitem{ainsworth2004dispersive}
Ainsworth, M.: Dispersive and dissipative behaviour of high order discontinuous
  {G}alerkin finite element methods.
\newblock Journal of Computational Physics \textbf{198}(1), 106--130 (2004)

\bibitem{AinsworthMonkMuniz}
Ainsworth, M., Monk, P., Muniz, W.: Dispersive and dissipative properties of
  discontinuous {G}alerkin finite element methods for the second-order wave
  equation.
\newblock J. Sci. Comput. \textbf{27}, 5--40 (2006)

\bibitem{DATH_UP}
Appel\"o, D., Hagstrom, T.: A new discontinuous {G}alerkin formulation for wave
  equations in second order form.
\newblock {SIAM} Journal On Numerical Analysis \textbf{53}(6), 2705--2726
  (2015)

\bibitem{appelo2018energy}
Appel{\"o}, D., Hagstrom, T.: An energy-based discontinuous {G}alerkin
  discretization of the elastic wave equation in second order form.
\newblock Computer Methods in Applied Mechanics and Engineering \textbf{338},
  362--391 (2018)

\bibitem{appelo2020energy}
Appel\"o, D., Hagstrom, T., Wang, Q., Zhang, L.: An energy-based discontinuous
  {G}alerkin method for semilinear wave equations.
\newblock Journal of Computational Physics \textbf{418}(109608) (2020)

\bibitem{Gdiff}
Banks, J., Hagstrom, T.: On {G}alerkin difference methods.
\newblock J. Comput. Phys. \textbf{313}, 310--327 (2016)

\bibitem{sipDGGD}
Banks, J.W., Buckner, B.B., Hagstrom, T., Juhnke, K.: Discontinuous {G}alerkin
  {G}alerkin differences for the wave equation in second-order form.
\newblock SIAM J. Sci. Comp.  (2021).
\newblock To appear

\bibitem{BANKS20125854}
Banks, J.W., Henshaw, W.D.: Upwind schemes for the wave equation in
  second-order form.
\newblock Journal of Computational Physics \textbf{231}(17), 5854--5889 (2012).
\newblock \doi{https://doi.org/10.1016/j.jcp.2012.05.012}.
\newblock
  \urlprefix\url{http://www.sciencedirect.com/science/article/pii/S0021999112002367}

\bibitem{chan2017weight}
Chan, J., Hewett, R.J., Warburton, T.: Weight-adjusted discontinuous {G}alerkin
  methods: curvilinear meshes.
\newblock SIAM Journal on Scientific Computing \textbf{39}(6), A2395--A2421
  (2017)

\bibitem{chan2017weight2}
Chan, J., Hewett, R.J., Warburton, T.: Weight-adjusted discontinuous {G}alerkin
  methods: wave propagation in heterogeneous media.
\newblock SIAM Journal on Scientific Computing \textbf{39}(6), A2935--A2961
  (2017)

\bibitem{ChouShuXing2014}
Chou, C.S., Shu, C.W., Xing, Y.: Optimal energy conserving local discontinuous
  {G}alerkin methods for second-order wave equation in heterogeneous media.
\newblock Journal of Computational Physics \textbf{272}, 88--107 (2014).
\newblock \doi{http://dx.doi.org/10.1016/j.jcp.2014.04.009}.
\newblock
  \urlprefix\url{http://www.sciencedirect.com/science/article/pii/S0021999114002721}

\bibitem{GSSwave}
Grote, M.J., Schneebeli, A., Sch{\"o}tzau, D.: Discontinuous {G}alerkin finite
  element method for the wave equation.
\newblock SIAM Journal on Numerical Analysis \textbf{44}(6), 2408--2431 (2006).
\newblock \urlprefix\url{http://www.jstor.org/stable/40232901}

\bibitem{GDDG}
Hagstrom, T., Banks, J.W., Buckner, B.B., Juhnke, K.: Discontinuous {G}alerkin
  difference methods for symmetric hyperbolic systems.
\newblock J. Sci. Comp. \textbf{81}, 1509--1526 (2019)

\bibitem{hesthaven2007nodal}
Hesthaven, J.S., Warburton, T.: Nodal discontinuous Galerkin methods:
  algorithms, analysis, and applications.
\newblock Springer Science \& Business Media (2007)

\bibitem{hu1999analysis}
Hu, F.Q., Hussaini, M., Rasetarinera, P.: An analysis of the discontinuous
  {G}alerkin method for wave propagation problems.
\newblock Journal of Computational Physics \textbf{151}(2), 921--946 (1999)

\bibitem{JonesPlassmannIC}
Jones, M., Plassmann, P.: Algorithm 740: {F}ortran subroutines to compute
  improved incomplete {C}holesky factorizations.
\newblock ACM Trans. Math. Soft. (TOMS) \textbf{21}, 5--17 (1995)

\bibitem{GDgeo}
Kozdon, J., Wilcox, L., Hagstrom, T., Banks, J.: Robust approaches to handling
  complex geometries with {G}alerkin difference methods.
\newblock J. Comput. Phys. \textbf{392}, 483--510 (2019)

\bibitem{lynch1964direct}
Lynch, R.E., Rice, J.R., Thomas, D.H.: Direct solution of partial difference
  equations by tensor product methods.
\newblock Numerische Mathematik \textbf{6}(1), 185--199 (1964)

\bibitem{Mattsson2012}
Mattsson, K.: Summation by parts operators for finite difference approximations
  of second--derivatives with variable coefficient.
\newblock J. Sci. Comput. \textbf{51}, 650--682 (2012)

\bibitem{Mattsson2004}
Mattsson, K., Nordstr\"{o}m, J.: Summation by parts operators for finite
  difference approximations of second derivatives.
\newblock J. Comput. Phys. \textbf{199}, 503--540 (2004)

\bibitem{moura2015linear}
Moura, R.C., Sherwin, S., Peir{\'o}, J.: Linear dispersion--diffusion analysis
  and its application to under-resolved turbulence simulations using
  discontinuous {G}alerkin spectral/hp methods.
\newblock Journal of Computational Physics \textbf{298}, 695--710 (2015)

\bibitem{IPDG_Elastic}
Riviere, B., Wheeler, M.: Discontinuous finite element methods for acoustic and
  elastic wave problems. part i: semidiscrete error estimates.
\newblock Contemporary Mathematics \textbf{329}, 271--282 (2003)

\bibitem{SchmitzYing3D}
Schmitz, P.G., Ying, L.: A fast nested dissection solver for {C}artesian 3{D}
  elliptic problems using hierarchical matrices.
\newblock J. Comput. Phys. \textbf{258}, 227--245 (2014)

\bibitem{SBPrev}
Sv\"ard, M., Nordstr\"om, J.: Review of summation-by-parts schemes for
  initial-boundary-value problems.
\newblock J. Comput. Phys. \textbf{268}, 17--38 (2014)

\bibitem{zhang2019energy}
Zhang, L., Hagstrom, T., Appel\"{o}, D.: An energy-based discontinuous
  {G}alerkin method for the wave equation with advection.
\newblock SIAM Journal on Numerical Analysis \textbf{57}(5), 2469--2492 (2019)

\end{thebibliography}
\bibliographystyle{spmpsci}

\end{document}